\documentclass[11pt]{article}
\usepackage[a4paper]{geometry}
\usepackage{amssymb}
\usepackage{amsmath,amsfonts}
\usepackage{mathtools}
\usepackage{amsthm}
\usepackage{latexsym}
\usepackage{mathrsfs}
\usepackage{color}
\usepackage{comment}
\usepackage{hyperref}
\usepackage{graphicx}
\usepackage{booktabs}
\usepackage{authblk}
\usepackage{subcaption}
\usepackage[ruled,vlined,linesnumbered,inoutnumbered]{algorithm2e} 
\usepackage{natbib}
 \bibpunct[, ]{(}{)}{,}{a}{}{,}%
\usepackage{bbm}
\usepackage{multirow}
\usepackage{bm}
\usepackage{cleveref}
\usepackage{xcolor}
\newcommand{\argmin}{\arg\min}
\newcommand{\argmax}{\arg\max}

\DeclarePairedDelimiter\ceil{\lceil}{\rceil}
\DeclarePairedDelimiter\floor{\lfloor}{\rfloor}
\newcommand{\Ascr}{\mathcal{A}}
\newcommand{\Bscr}{\mathcal{B}}
\newcommand{\Cscr}{\mathcal{C}}
\newcommand{\Dscr}{\mathcal{D}}

\newcommand{\Fscr}{\mathcal{F}}
\newcommand{\Gscr}{\mathcal{G}}

\newcommand{\Lscr}{\mathcal{L}}

\newcommand{\Nscr}{\mathcal{N}}

\newcommand{\Sscr}{\mathcal{S}}

\newcommand{\Xscr}{\mathcal{X}}

\newcommand{\mE}{\mathbb{E}}

\newcommand{\mR}{\mathbb{R}}

\newcommand{\mZ}{\mathbb{Z}}

\begin{document}

\title{Artificial Intelligence for Operations Research: Revolutionizing the Operations Research Process}


\author{Zhenan Fan, Bissan Ghaddar, Xinglu Wang, Linzi Xing, Yong Zhang, Zirui Zhou \thanks{The authors are listed in alphabetical order}}


\date{}
\maketitle
\begin{abstract}
    The rapid advancement of artificial intelligence (AI) techniques has opened up new opportunities to revolutionize various fields, including operations research (OR). This survey paper explores the integration of AI within the OR process (AI4OR) to enhance its effectiveness and efficiency across multiple stages, such as parameter generation, model formulation, and model optimization. By providing a comprehensive overview of the state-of-the-art and examining the potential of AI to transform OR, this paper aims to inspire further research and innovation in the development of AI-enhanced OR methods and tools. The synergy between AI and OR is poised to drive significant advancements and novel solutions in a multitude of domains, ultimately leading to more effective and efficient decision-making. 
\end{abstract}
\maketitle

\section{Introduction}
Operations Research (OR) is an interdisciplinary field that employs advanced analytical techniques and methodologies to support decision-making processes in organizations, aiming to improve efficiency, optimize resource allocation, and achieve desired objectives. By leveraging mathematical models, optimization algorithms, simulation, and statistical methods, OR aids in addressing complex problems in various domains, including logistics, supply chain management, transportation, energy, manufacturing, finance, healthcare, and public services, among others \citep{Winston2004}. 
The general framework for operations research includes the following steps~\citep{rajgopal2004principles}
\begin{itemize}

\item Problem identification and definition: The initial step requires a thorough understanding of the problem, including its context, objectives, and relevant constraints. This involves engaging with stakeholders to clarify their expectations, goals, and requirements, as well as identifying any trade-offs or conflicting priorities that may arise during the decision-making process. A well-defined problem statement provides a solid foundation for the subsequent stages of the OR process.

\item Parameter Generation: The next step involves generating key parameters in the optimization model, such as the objective coefficients and the constraint matrix. In practice, we might have relevant data from various sources, such as historical records, expert opinions, market research, or sensor readings. The collected data may need preprocessing, such as cleaning, normalization, aggregation, or transformation, to ensure its quality and suitability for the subsequent modeling stage. These data are then converted to modeling parameters by human experts or AI models. For example, in a supply chain planning problem, although the expert has not yet written down explicitly the objective and constraint equations, they can first decide on key parameters such as supply, demand, and profit according to their understanding of the problem and analysis of data. 

\item Model formulation: With a clear understanding of the problem and the corresponding key parameters, a mathematical or simulation model is developed to represent the system under investigation. Models are simplifications of reality, designed to capture the essential elements of the problem while abstracting away unnecessary details. Depending on the problem's nature and requirements, various modeling techniques can be employed, including linear/non-linear programming, integer programming, stochastic programming, queuing theory, network models, game theory, or agent-based simulations, among others. Before applying the model in practice, its validity must be assessed to ensure that it accurately captures the problem's essential characteristics and adheres to the system's logical and physical constraints.

\item Model optimization: Once the model has been validated, it is solved and analyzed using various techniques to identify optimal or near-optimal solutions that satisfy the problem's objectives and constraints. This may involve employing optimization algorithms, heuristics, metaheuristics, or simulation-based methods, depending on the model's complexity and the desired level of solution quality.

\item Interpretation and validation: In the final stage of the OR process, the implemented solution is reviewed and evaluated to determine its effectiveness in meeting the desired objectives. If necessary, the model may be updated, and the process iterated to further refine the solution and enhance the system's performance. This step embodies the iterative and ongoing nature of operations research as a tool for continuous improvement.

\end{itemize}

While these steps have proven effective in the past, recent advances in artificial intelligence (AI) are poised to revolutionize the way we approach and solve OR problems. AI techniques have the potential to enhance every stage of the OR process, facilitating the development of more accurate and efficient models and offering innovative solutions to complex problems. In this survey paper, we explore AI4OR, i.e., how AI can be combined with OR, with a focus on three key aspects: parameter generation, model formulation, and model optimization.

First, in the parameter generation phase, AI can be employed to improve the quality and relevance of data used to formulate the mathematical models. One promising approach is the predict-then-optimize framework of \cite{elmachtoub2022smart}, which leverages AI techniques to make data-driven predictions about uncertain variables in the decision-making process. By incorporating advanced AI algorithms, such as deep learning and reinforcement learning, predict-then-optimize can efficiently handle high-dimensional and complex data structures, as well as adapt to dynamic environments. Furthermore, AI-powered feature selection and dimensionality reduction techniques can be utilized to identify the most important variables and relationships within the data, leading to more parsimonious and interpretable models.

Second, in the model formulation phase, AI can be used to bridge the gap between natural language problem description and mathematical models. This is particularly relevant in situations where domain experts may struggle to translate their knowledge into mathematical terms. Large language models, such as ChatGPT \citep{brown2020language} and Llama \citep{touvron2023llama}, have shown remarkable success in understanding and generating natural language and can be employed to automatically convert problem descriptions into mathematical formulations~\citep{bubeck2023sparks,ramamonjison2023nl4opt}. By harnessing the power of natural language processing (NLP) and AI, these models can extract relevant information, identify key constraints and objectives, and translate them into a mathematical representation suitable for optimization.

Third, in the model optimization phase, AI techniques can be exploited to enhance the performance of optimization algorithms, shifting away from traditional methods towards a more adaptive, learning-driven approach. Classic optimization methods, such as gradient descent, conjugate gradient, Newton steps, and branch-and-bound, are constructed based on theoretical foundations and the implementation of optimization experts. While these methods offer performance guarantees, they may not always be the most efficient or effective solutions for specific problem instances. This survey paper focuses on three main categories in this direction: automatic algorithm configuration, continuous optimization algorithm selection and design, and discrete optimization algorithm selection and design. A brief description of these three categories is presented next:

\begin{enumerate}
    \item Automatic Algorithm Configuration: This category encompasses techniques that use AI to fine-tune the parameters of existing optimization algorithms, enabling better performance on specific problem instances. Methods such as Bayesian optimization, genetic algorithms, and reinforcement learning can be employed to intelligently search the parameter space and identify configurations that yield improved performance \citep{ansotegui2009gender,JMLR:v23:21-0888,AnaHoo20}.
    \item Continuous Optimization Algorithm Selection and Design: 
    AI techniques are being employed to enhance the optimization algorithm that solves problems with continuous variables. Techniques like learning to optimize~\citep{wichrowska2017learned}, the adaptive penalty for ADMM~\citep{zeng2022reinforcement}, and smart column selection~\citep{chi2022deep} can be utilized to dynamically determine the most appropriate step size, balance exploration and exploitation trade-offs, and accelerate the optimization process.
    \item Discrete Optimization Algorithm Selection and Design: This category focuses on the application of AI techniques in optimization problems with discrete variables, such as those encountered in combinatorial optimization problems. AI-driven heuristics~\citep{di2016dash}, metaheuristics~\citep{talbi2009metaheuristics}, and learning-based approaches~\citep{gomory1960algorithm} can be employed to enhance algorithms like branch-and-bound and cutting-plane methods, and improve solving mixed integer programming problems. 
\end{enumerate}

Our survey covers important stages in the pipeline of OR and investigates how AI can assist each stage of the pipeline, which offers a holistic perspective and allows us to gain valuable insights into the potential synergies between AI and OR. Our key contributions involve two aspects. First, we comprehensively examine the different components of the OR pipeline while existing surveys focus solely on specific aspects. For instance, 
\citet{bengio2021machine} focus on the integration of AI with combinatorial optimization problems, enabling autonomous learning and decision-making on a chosen specific set of problems. 
\citet{Zhang2022ASF} review how AI assists Mixed Integer Programming (MIP) algorithms including branch-and-bound and heuristic methods. 
\citet{Lodi2017OnLA} focus on the AI-enhanced variable and node selection in the branch-and-bound algorithm for MIPs. 
\citet{schede2022survey} survey automated algorithm configuration methods. 
\citet{Kotary2021EndtoEndCO} survey two directions that leverage AI for constrained optimization (CO): AI-augmented CO, which enhances optimization algorithms with AI assistance, and End-to-End CO learning, where machine learning directly predicts the solution of CO. Additionally, some sections in our survey have not been discussed before in the existing surveys. For example, the model formulation (see Section \ref{sec:model_formulation}) and enhancement (see Section \ref{sec:admm} and \ref{sec:column_generation}) for specific algorithms like ADMM and column generation. 

Second, we analyze the pipeline as a whole and emphasize the interactions between different components. 
For example, the interaction between mathematical model parameter generation and optimization attracted increased attention from recent literature. 
Compared with the traditional predict-then-optimize paradigm that isolates the two components, recent approaches started to explore their interactions.  
The smart predict-then-optimize paradigm by \cite{elmachtoub2022smart} and \cite{ amos2017optnet} gathers feedback from later decision errors to refine the prediction of parameters. 
The ``integrated prediction and optimization" paradigm presented by \cite{bertsimas2020from, mip-constraint} and \cite{bergman2020janos} presents another possible interaction. Recall that the traditional predict-then-optimize paradigm involves first predicting the mathematical model parameters and then deriving an optimal decision. 
Thus, the parameter is not dependent on the decision. 
In contrast, the new paradigm allows the parameters to be affected by the decision, i.e., the parameter generation takes the future decision into account. 

We survey this emerging direction, describe existing works (See Section \ref{sec:param_gen}), and envision other potential interactions between different components within the optimization pipeline (See Section \ref{sec:conclusion}). 
It is essential to clarify that our survey's scope is limited to methods that involve the use of optimization software. 
There is a subset of real-world operational problems that do not necessitate optimization software, e.g., unconstraint problem, time series prediction, or classification. 
For these problems, end-to-end AI solutions is possible \citep{kraus2020deep,Kotary2021EndtoEndCO,Zhang2023ARO,Yan2020LearningFG,Guo2019SolvingCP}. 
However, this survey does not cover these problems and their end-to-end AI solutions. 
\citet{kraus2020deep} review the before-mentioned problems, like predicting the movements of stocks, from an operational point of view and how AI will achieve high prediction performance. 
Considering that AI solutions are often a black box, \citet{debock2023explainable} disucss how the integration of explainability and ethical consideration can be taken into account in AI solutions alongside performance. 
Furthermore, while our focus is on how AI can boost operational research, it is well known that machine learning itself roots deeply in mathematical optimization. \citet{gambella2021optimization} survey mathematical optimization models presented in various AI algorithms, such as classification, clustering, deep learning, and Bayesian network structure learning.

\begin{figure}[th]
    \centering
    \includegraphics[width=1.01\linewidth]{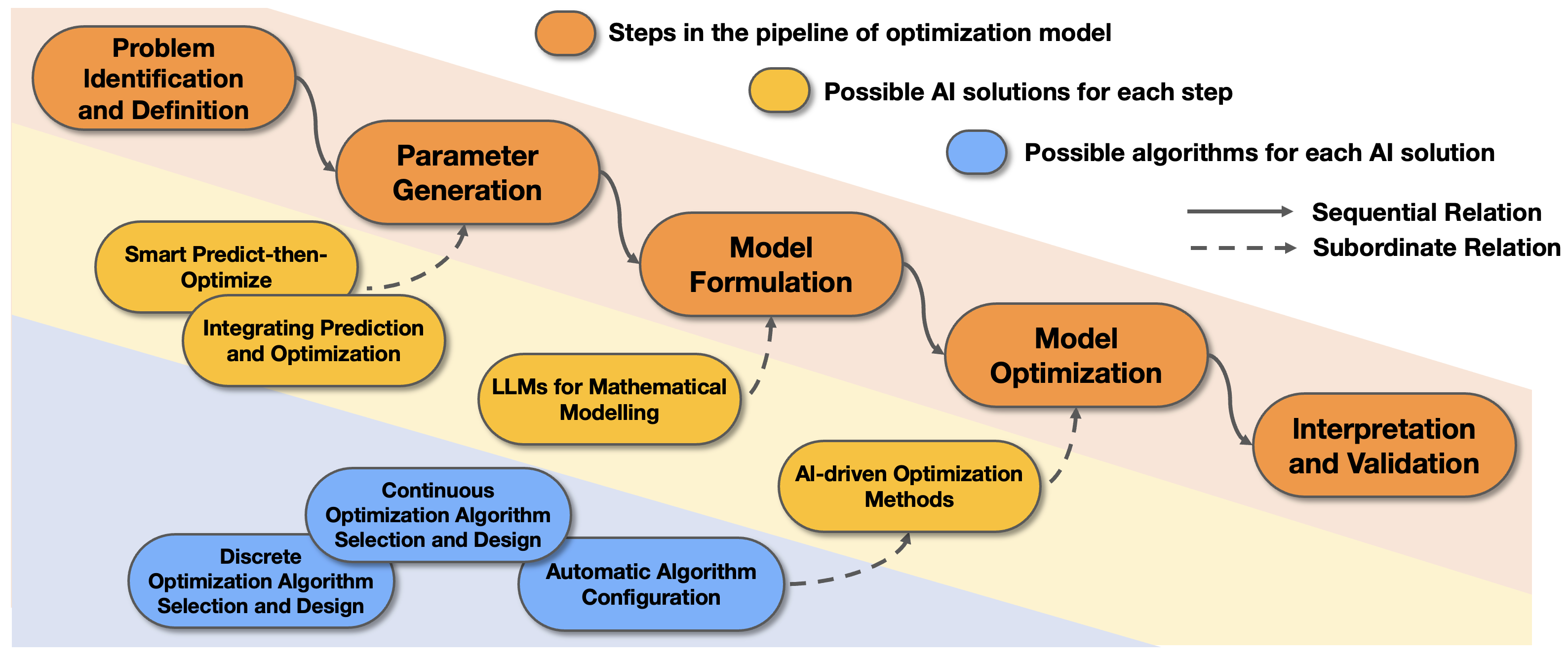}
    \caption{AI techniques in the OR process}
    \label{fig:ML4Solver}
\end{figure}

In summary, 
operation research has benefited from the advancement in various algorithms (see Table~\ref{tab:trad-algos}), computing power as well as commercial and open-source solvers. 
Solvers like \cite{gurobi}, \cite{cplex2009v12}, and OptVerse~\cite{Huawei2021} are generally applicable to various real-world applications and can solve large-scale problems efficiently. 
Meanwhile, the advances in AI led to a groundbreaking shift in the development of optimization algorithms. Instead of relying on structured algorithmic development, AI techniques learn from data (e.g., past solving experience), enhance existing methods or even create entirely novel solution methods. 
For specific classes of optimization problems, 
integration of artificial intelligence into operations research can even achieve better performance than existing solvers \citep{song2020general,Gasse2019,Khalil2016,jung2022learning,hutter2010automated,hutter2009paramils}. 
AI4OR promises to bring about a new era of innovation and efficiency. This survey paper will delve into the various AI techniques that can be employed at each stage of the OR process, providing a comprehensive overview of the state-of-the-art and exploring the potential of AI to transform the way we approach and solve complex decision-making problems. As we continue to develop and refine AI technologies, the synergy between AI and OR will undoubtedly lead to exciting advancements and novel solution methods in a multitude of domains. An overview illustrating how AI techniques help each step in the OR pipeline is shown in Figure~\ref{fig:ML4Solver}. 

\begin{table}[th]
    \centering
    \begin{tabular}{|c|c|c|}
    \hline
      Algorithm                         &  Target problem     & Assumption\\ \hline
      Gradient-based methods            &   $\min\limits_{x \in\mR^n}\enspace f(x)$                                                                     & $f$ is differentiable \\ 
      ADMM-type methods                 &   $\min\limits_{x \in \mR^n, y \in \mR^m}\enspace f(x) + g(s) \enspace\text{s.t.}\enspace Ax + s = b$         & $f$ and $g$ are continuous \\ 
      Simplex method                    &   $\min\limits_{x \in \mR^n}\enspace c^Tx \enspace\text{s.t.}\enspace Ax = b, \ x \geq 0$                     & $A$ has full rank\\
      Column-generation method          &   $\min\limits_{x\in\mR^{|\Omega|}} \enspace \sum\limits_{p\in\Omega}c_px_p \enspace\text{s.t.}\enspace \sum\limits_{p\in\Omega}x_p a_p = b, \ x \geq 0$                     & $|\Omega|$ is very large\\
      Cutting-plane method              &   $\min\limits_{x \in \mZ^n}\enspace c^Tx \enspace\text{s.t.}\enspace Ax = b, \ x \geq 0$                     & $A$ has full rank\\
      Branch-and-bound method           &   $\min\limits_{x \in \Xscr}\enspace f(x)$                                                                    & $\Xscr$ is finite\\
    \hline
    \end{tabular}
    \caption{Summary of involved algorithms and their target optimization problems}
    \label{tab:trad-algos}
\end{table}

\section{Preliminary: AI techniques}

There are two main challenges in operation research, and AI techniques have the potential to address them.  We start by specifically introducing these two challenges:

{1)} {The complex interaction among decision variables and constraints \citep{Cappart2021,Steever2019AnIA}:} Simple linear algebra-based heuristics may be able to identify that the constraint matrix has a block structure, then the original problem can be decomposed into several subproblems, and thus accelerate solving the optimization problem. However, when the interaction within the constraint matrix is more complex, e.g., more variables or constraints are coupling with each other, simple heuristics are either not applicable or worsen the computational performance. 
It is common that a group of optimization problems share certain characteristics that cannot be described (or easily obtained) mathematically. 
In such situations, AI tools are more applicable and efficient. 

{2)} {The computation cost of solving complex optimization problems \citep{krentel1986complexity}:} Many interesting and practical optimization problems are NP-hard. Even for polynomial-time solvable problems, e.g., LP and QP, when the problems scale up, the computational time grows quickly. 
Thanks to advances in hardware, software, and optimization solvers, many optimization problems can now be solved within acceptable time frames. 
However, due to the increase in data availability optimization problems often continue to grow in terms of scale and complexity that surpass the capabilities of existing solvers for real-world applications. 
Thus we need to harness AI techniques to further accelerate and improve the efficiency of solving OR models. 

In the following, we present brief intuitions about some of the AI techniques frequently used in OR and why they are helpful in addressing the challenges.
We will focus on two key aspects of AI techniques: 1) the {AI models} themselves and 2) the {learning algorithms} of these models. 

1) Commonly used {AI models} include  Graph Neural Networks (GNN,  Section \ref{sec:gnn}) and Recurrent Neural Network (RNN, Section \ref{sec:rnn}). 
GNN excels at handling the complex interactions in graphs, such as the real-world graphs represented in the parameter generation stage  \citep{elmachtoub2022smart} and the equivalent graph representation of LP, QP, or MILP problems in the AI-driven optimization stage. 
RNN is capable of retaining information from previous time steps. When applying RNN to an iterative algorithm, a time step corresponds to an iteration. 
Many iterative optimization algorithms are typically slow because they need to check one solution in each iteration and determine the subsequent solutions to assess. An inefficient or suboptimal choice can hinder the optimization progress. 
Given RNN's proficiency with sequential information (e.g., previous algorithm decisions and states), it is a helpful tool in such situations. 

2) Prominent {learning algorithms} include Reinforcement learning (see Section~\ref{sec:rl}) and Imitation Learning (see Section \ref{sec:imtate}). 
Reinforcement learning offers a solution to the challenge of high computational cost. 
As previously mentioned, this cost often stems from an inefficient or suboptimal decision made in early iterations. 
To tackle this, two more challenges appear: 
i) The metrics, like execution time or duality gap, are non-differentiable w.r.t the decisions.
ii) The need for a learning algorithm that permits long-term rewards to influence earlier decisions, termed as delayed rewards. 
Reinforcement learning effectively meets both these criteria. 
Imitation Learning is another strategy for the challenge of high computational cost. 
By design, imitation learning emulates expert behaviour. In optimization fields, these "experts" often represent computationally intensive methods well-documented in the literature. Instead of traditional input-label pairs as in supervised learning, imitation learning utilizes state-action pairs derived from expert demonstrations. The learning objective becomes the alignment of the model's predictions with the expert behaviour, which is differentiable. 
Upon completion, imitation learning yields a model that rapidly predicts expert behaviour, thereby effectively tackling the computational challenge.

In the subsequent portion of this section, we will give a detailed overview of these four AI techniques.  

\subsection{Graph Neural Networks} \label{sec:gnn}

\begin{figure}[t]
    \centering
    \includegraphics[width=.7\linewidth]{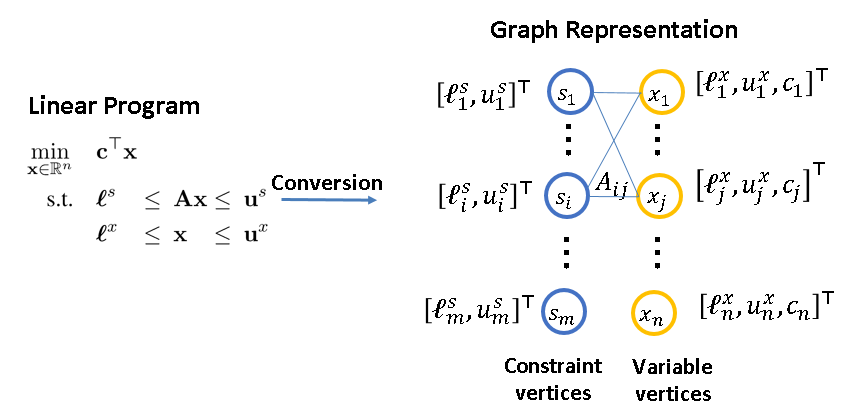}
    \caption{Graph representation for linear program. }
    \label{fig:lp-graph}
\end{figure}

\begin{figure}[b]
    \centering
    \includegraphics[width=.7\linewidth]{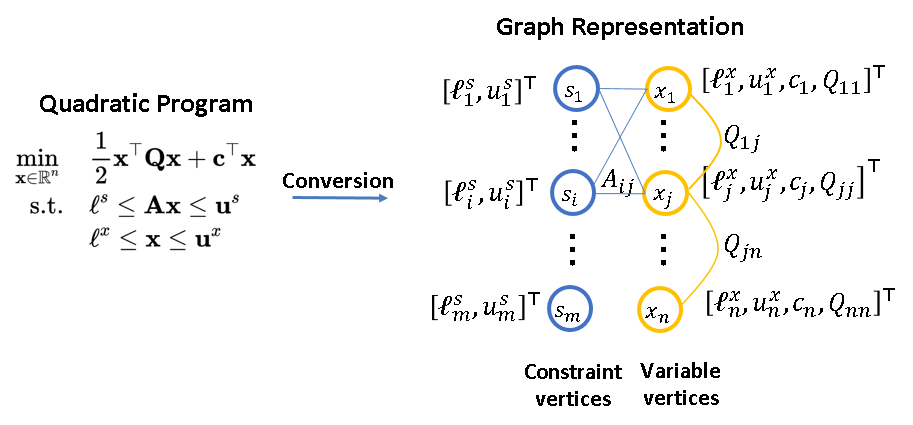}
    \caption{Graph representation for quadratic programs.}
    \label{fig:qp-graph} 
\end{figure}

Graph Neural Networks (GNNs) are a class of neural networks designed for processing graph-structured data \citep{zhou2020graph,wu2020comprehensive,xu2018powerful,velivckovic2023everything} . Consider a graph $G = (V, E)$, where $V = \{v_1, \dots, v_N\}$ is the set of nodes and $E \subseteq V \times V$ is the set of edges. Each node $v_i$ is associated with a feature vector $\bm{x}_i \in \mathbb{R}^d$. GNN aims to learn a mapping that maps the graph $G$ and its node features to a set of output vectors $\bm{y}_1, \dots, \bm{y}_N$, where $\bm{y}_i \in \mathbb{R}^m$. GNN typically consists of multiple layers, and the node representations are updated at each layer. Let $\bm{h}_i^{(l)}$ be the hidden representation of node $v_i$ at the $l$-th layer of the GNN. The initial node representations are given by the input features: $\bm{h}_i^{(0)} = \bm{x}_i$. The GNN updates the node representations in each layer as follows:
\begin{equation*}
  \bm{h}_i^{(l+1)} = \phi^{(l)} \left( \bm{h}_i^{(l)}, \sum_{v_j \in \mathcal{N}(v_i)} \psi^{(l)}(\bm{h}_i^{(l)}, \bm{h}_j^{(l)}) \right),
\end{equation*}
where $\phi^{(l)}$ and $\psi^{(l)}$ are neural network functions with learnable parameters, and $\mathcal{N}(v_i)$ is the set of neighbors of node $v_i$ in the graph. The updated node representations are then used to compute the output vectors:
\begin{equation*}
  \bm{y}_i = \omega(\bm{h}_i^{(L)}),
\end{equation*}
where $L$ is the number of layers in the GNN, and $\omega$ is a neural network function with learnable parameters. The parameters of the GNN are given by $\bm{\theta} = \{\phi_\theta^{(l)}, \cdots, \theta^{(l)}, \cdots, \phi_\theta^{(L)}, \omega_\theta \}$, where $\theta^{(l)} = \{ \phi_\theta^{(l)}, \psi_\theta^{(l)} \}$ is the parapemter for $l$-th layer. 

In summary, GNN processes graph-structured data by updating the node representations through multiple layers using neural network functions. The GNN is defined by a set of learnable parameters and is trained by minimizing a loss function that measures the difference between the predicted output vectors and the ground truth output vectors using an optimization algorithm. As shown in Figures~\ref{fig:lp-graph} and \ref{fig:qp-graph} (these figures are adapted and simplified from \cite{jung2022learning,fan2023smart} and \cite{Gasse2019}), 
optimization problems like linear programming and quadratic programming can be equivalently transformed into a graph. 
This transformation is beneficial since the complex interactions between constraints and variables are expressed, and the permutation invariance property is preserved \citep{velivckovic2023everything,wu2020comprehensive}. 
Specifically, in the context of optimization, permutation invariance implies that rearranging the order of constraints or variables does not affect the problem's equivalence or the solution's validity. This property is retained in the graph representation, i.e., permuting the constraints and variables results in an equivalent graph and does not impact the prediction made by GNN.
For representing mixed integer programming, we can simply add another attribute to the variable nodes denoting whether the variables are integer or continuous.  More complex features like the sparsity of columns can also be added \citep{fan2023smart}.  Thus, GNN is useful for representing optimization problems and providing feature representation when we want to accelerate the model optimization process with AI techniques. 

\subsection{Recurrent Neural Networks} \label{sec:rnn}
Recurrent Neural Networks (RNNs) are a class of neural networks designed for processing sequential data \citep{medsker2001recurrent,yu2019review,staudemeyer2019understanding}. 
Let $x^1, \dots, x^T$ be a sequence of input vectors, where $T$ is the length of the sequential data and $x^t \in \mathbb{R}^d$ represents the input vector at $t$ step in the sequence. 
The RNN consists of neural network functions $f(\theta_f; \cdot)$ and $g(\theta_g;\cdot)$, such that we have 
\begin{equation}
\begin{aligned}
    h(t) &= f(\theta_f; x^t, h^{t-1}) \\
    y(t) &= g(\theta_g; h^t)
\end{aligned} \label{eq:rnn}
\end{equation}
where parameter $\bm{\theta}=\{\theta_f, \theta_g\}$ are the learnable parameters and $h^t$ is the hidden state at time step $t$. 

For standard RNN, $f$ and $g$ can be Multilayer Perceptron (MLP) models or even simply a single matrix multiplication operation. 
However, standard RNNs face a significant challenge with the vanishing gradient problem, which makes them struggle to capture long-term dependencies in the sequential data. 
To illustrate, consider the case that function $f$ is a matrix multiplication operation. The initial input information $x^0$ is subjected to repeated multiplication with the same matrix $\theta_f$ over successive time steps.
The cumulative effect of this iterative multiplication can cause the contribution of $x^0$ towards $h^t$ decreasing geometrically, depending on the eigenvalues of the matrix $\theta_f$. 
This reduction of influence shows that RNN could `forget' distant past information. 
To overcome this limitation, Long Short-Term Memory Networks (LSTMs)  were introduced \citep{yu2019review,staudemeyer2019understanding}. LSTMs are a special kind of RNN. The functions $f$ and $g$ are structured differently to include a gating mechanism. This gating mechanism is designed to control and manage the flow of information within the neural network, enabling LSTMs to retain information for much longer periods of time and to better capture long-term dependencies.

Ideally, the sequence $y^1, \dots, y^T$ generated by the RNN with parameter $\bm{\theta}$ should match the ground-truth sequence $z^1, \dots, z^T$ annotated by human experts. 
The exact match is challenging and unnecessary, but it motivates us to adopt the difference between these two sequences as the loss function $\mathcal{L}(\bm{\theta})$ for training an RNN model. 
\begin{equation*}
\mathcal{L}(\bm{\theta}) = \frac{1}{T} \sum_{t=1}^T \ell(y^t, z^t), \text{ where } y^t \text{ depends on } \bm{\theta}
\end{equation*}
$\ell(\cdot, \cdot)$ is some user-specified function measuring the difference between $y^t$ and $z^t$. The RNN is trained by minimizing the loss function $\mathcal{L}(\bm{\theta})$ with respect to the learnable parameters $\theta$ using an optimization algorithm, such as stochastic gradient descent or one of its variants. 
Since $\theta_f$ and $\theta_g$ are reused at each step $t$ in Equation~\eqref{eq:rnn}, 
The optimization process typically involves backpropagation through time \citep{medsker2001recurrent}, which is an extension of the standard backpropagation algorithm to handle recurrent neural networks.

In summary, an RNN consists of a sequence of input vectors, hidden states, and output vectors. The network is defined by a set of learnable functions. The RNN is trained by minimizing a loss function that measures the difference between the predicted output sequence and the ground truth output sequence using an optimization algorithm. 
Due to RNN's ability to operate on sequential information, it is useful for accelerating iterative optimization algorithms, where historical states and decisions matter for the current decision.

\subsection{Reinforcement Learning} \label{sec:rl}
Reinforcement learning is an AI technique in which an agent learns to make decisions through interactions with its environment. The agent seeks to maximize the cumulative reward it receives over time. The problem is typically formulated as a Markov Decision Process (MDP), defined by the tuple $(S, A, P, R, \gamma)$, where:

\begin{itemize}
    \item $S=\{s_1,\cdots, s_n\}$ is the sequence of states.
    \item $A=\{a_1,\cdots, a_n\}$ is the sequence of actions.
    \item $P(s_{t+1} | s_t, a_t)$ is the state transition probability function, describing the probability of transitioning from state $s_t$ to state $s_{t+1}$ after taking action $a_t$.
    \item $R(s_t, a_t, s_{t+1})$ is the reward function, providing the immediate reward received by the agent after taking action $a_t$ in state $s_t$ and transitioning to state $s_{t+1}$.
    \item $\gamma \in [0, 1]$ is the discount factor, which determines the importance of future rewards relative to immediate rewards.
\end{itemize}

A policy, denoted as $\pi(a|s)$, is a probability distribution over actions given the current state. The goal of reinforcement learning is to find an optimal policy $\pi^*$ that maximizes the expected return from any initial state.
To understand this objective, we first introduce the {discounted return} at the time step \( t \):
\[
G(s_t,a_t) = \sum_{k=0}^{\infty} \gamma^k R(s_{t+k}, a_{t+k}, s_{t+k+1})
\]
Given the stochastic nature of state transitions, starting from an initial state \( s_t \) can lead to different possible future states. 
In other words, \( \{ s_{t+1}, s_{t+2}, \ldots \} \) and $G(s_t,a_t)$ are all random variables. 
Therefore, the goal is actually to maximize the expected value of the discounted return from any state $s$, commonly referred to as {expected return}, which can be formally defined as:
\begin{equation} \label{eq:exp-return}
    V^\pi(s) = \mathbb{E}_{\pi} \left[ G(s_0,a_0) \big| s_0 = s \right]
\end{equation}
Here, $a_0$ is also a random variable decided by $\pi(a_0|s_0)$.
This expected return from state $s$ measures the quality of state $s$. Thus, $V^\pi(s) $ is also named as the state value function. 


To find the optimal policy, there are two streams of methods: 1) value function-based methods and 2) policy optimization. In the value function-based method, 
a commonly used tool is the state-action value function, denoted as \( Q^\pi(s, a) \). Formally, this function is defined as follows:
\[
Q^\pi(s, a) = \mathbb{E}_{\pi} \left[ G(s_0, a_0) \big| s_0 = s, a_0 = a \right]
\]
This function serves to quantitatively evaluate the quality of taking action \( a \) in state \( s \). 
Given \( Q^\pi(s, a) \), the optimal policy to choose the action given by \( \argmax_a Q^\pi(s, a) \).



Policy optimization aims to find the optimal policy $\pi^*$ that maximizes the expected return 
given in Equation \eqref{eq:exp-return}. 
One common approach is to parameterize the policy using a function approximator, such as a neural network, with learnable parameters $\theta$. The policy is then represented as $\pi_\theta(a|s)$. 
Now the expected return objective $V^{\pi_\theta}(s) $ is parameterized by $\theta$. 
For simplicity, we denote this objective as $J(\theta)=V^{\pi_\theta}(s)$.

The policy gradient method is a popular approach to optimize the policy parameters $\theta$. It computes the gradient of the objective function with respect to the policy parameters and updates the parameters using gradient ascent. The policy gradient can be expressed as:
\begin{equation} 
\nabla_\theta J(\theta) = \mathbb{E}_{\pi_\theta} \left[ \sum_{t=0}^T \nabla_\theta \log \pi_\theta(a_t|s_t) G(s_t,a_t) \right].
\end{equation}


The policy parameters $\theta$ are updated using a gradient-based optimization algorithm, such as stochastic gradient ascent or an adaptive optimization method like Adam (\citep{ruder2016overview,kingma2014adam}:
\[
\theta \leftarrow \theta + \alpha \nabla_\theta J(\theta),
\]
where $\alpha$ is the learning rate.

A number of improving techniques have been developed to enhance the basic policy gradient method. Their purpose is to stabilize the learning process, improve convergence, and boost overall performance. These are a few popular approaches:
\begin{itemize}
\item \textit{Trust region}: 
The principle is to maintain a ``trust region" for the policy that the agent outputs \citep{TRPO}. In essence, the agent creates a new policy that isn't too deviant from the old policy, confining the updates within a predefined boundary or ``trust region." This technique is commonly used when there is a need for stability in the learning process is crucial.
\item \textit{Surrogate objective}: This strategy substitutes the original objective function with a surrogate function, known as PPO algorithm \citep{PPO}. This method is less complex to implement compared to the trust region approach, and it also promotes stable training.
\item \textit{Value function approximations}: In this strategy, value function approximations are incorporated as a baseline to enhance policy optimization \citep{A2C, A3C, SAC}. This method capitalizes on both the learnable value function and policy optimization process, usually resulting in improved empirical performance.
\end{itemize}

In conclusion, reinforcement learning equips an agent with a reward that doesn't necessarily need to be differentiable regarding the agent's decisions. Moreover, the agent will not be short-sighted since the learning objective is maximizing {cumulative} rewards. 
Given the continuous advancements in the state of the art, reinforcement learning is becoming increasingly stable and applicable to real-world scenarios. 

\subsection{Imitation learning} \label{sec:imtate}
Imitation learning is a type of learning algorithm where an agent learns to make decisions by observing demonstrations provided by an expert \citep{hussein2017imitation}. The goal is to learn a policy that mimics the expert's behaviour as closely as possible. 
This approach is particularly useful in scenarios where reinforcement learning methods face challenges, such as when designing an appropriate reward function is difficult, when environmental interactions are costly, or when the learning process requires a large number of samples. 
There are several methods for imitation learning, including behavior cloning by \cite{codevilla2019exploring} and \cite{torabi2018behavioral} and interactive imitation learning by \citep{forward-train}. {Behavior cloning} is a simple and widely used imitation learning approach \citep{codevilla2019exploring,torabi2018behavioral}. Given a dataset of expert demonstrations $\mathcal{D} = \{(s_i, a_i)\}_{i=1}^N$, where $s_i$ is a state and $a_i$ is the corresponding expert action, the goal is to learn a policy $\pi_\theta(a|s)$ that maps state to actions. The policy is usually represented by a function approximator, such as a neural network, with learnable parameters $\theta$. The learning problem is formulated as a supervised learning task, minimizing the loss function $\mathcal{L}(\theta)$:

\[
\mathcal{L}(\theta) = \frac{1}{N} \sum_{i=1}^N \ell(a_i, \pi_\theta(a|s_i)),
\]
where $\ell$ is a distance or divergence measure between the expert action $a_i$ and the predicted action $\pi_\theta(a|s_i)$. Common choices for $\ell$ include the mean squared error for continuous actions or the cross-entropy loss for discrete actions.

Despite its simplicity, behaviour cloning can suffer from ``cascading errors" \citep{forward-train}. This happens because the model's training data consists of sequences of actions and states from the expert. If, at any point, the model deviates from the expert's actions, it may find itself in a state that it has never seen during training. This can lead to incorrect actions in the next step, which in turn leads to more unfamiliar states, causing the errors to cascade. 
To address this issue, \citet{forward-train} propose {interactive imitation learning}. After the model takes a new action, the pair of the current state and expert action (if it exists) is added to demonstrations $\Dscr$. 
Despite the annotation burden on the expert increases, 
the model learns to handle the states deviating from the original expert's actions, helping to correct its mistakes and prevent cascading errors. 
Imitation learning is closely related to reinforcement learning, as both aim to learn policies that maximize some objective. 
While RL seeks to maximize the cumulative reward directly, imitation learning aims to learn from expert demonstrations. 

Imitation learning is often considered more ``sample efficient," {i.e.}, it learns effectively from a smaller number of expert demonstrations. This is because imitation learning leverages the knowledge of an expert, bypassing the need for the extensive ``trial-and-error" process that characterizes RL. This trial-and-error process involves the agent making numerous attempts at a task, learning from its mistakes and successes, which can be time-consuming and resource-intensive. However, the performance of imitation learning is bounded by the capabilities of the expert. For instance, if the expert prioritizes short-term benefits, the agent trained through imitation learning may also be short-sighted. 
In such cases, RL is better suited because it allows the agent to learn from its own experiences, rather than strictly following an expert. This means that the agent can adapt based on the consequences of its own actions, taking into account both immediate and future consequences.

In conclusion, the presented techniques have specific benefits for operational research. Each algorithm has its own suitable application scenario, as we will see in the coming sections.

\section{Model Paramater Generation} \label{sec:param_gen} 
 
\begin{figure}
    \centering
    \includegraphics[width=\linewidth]{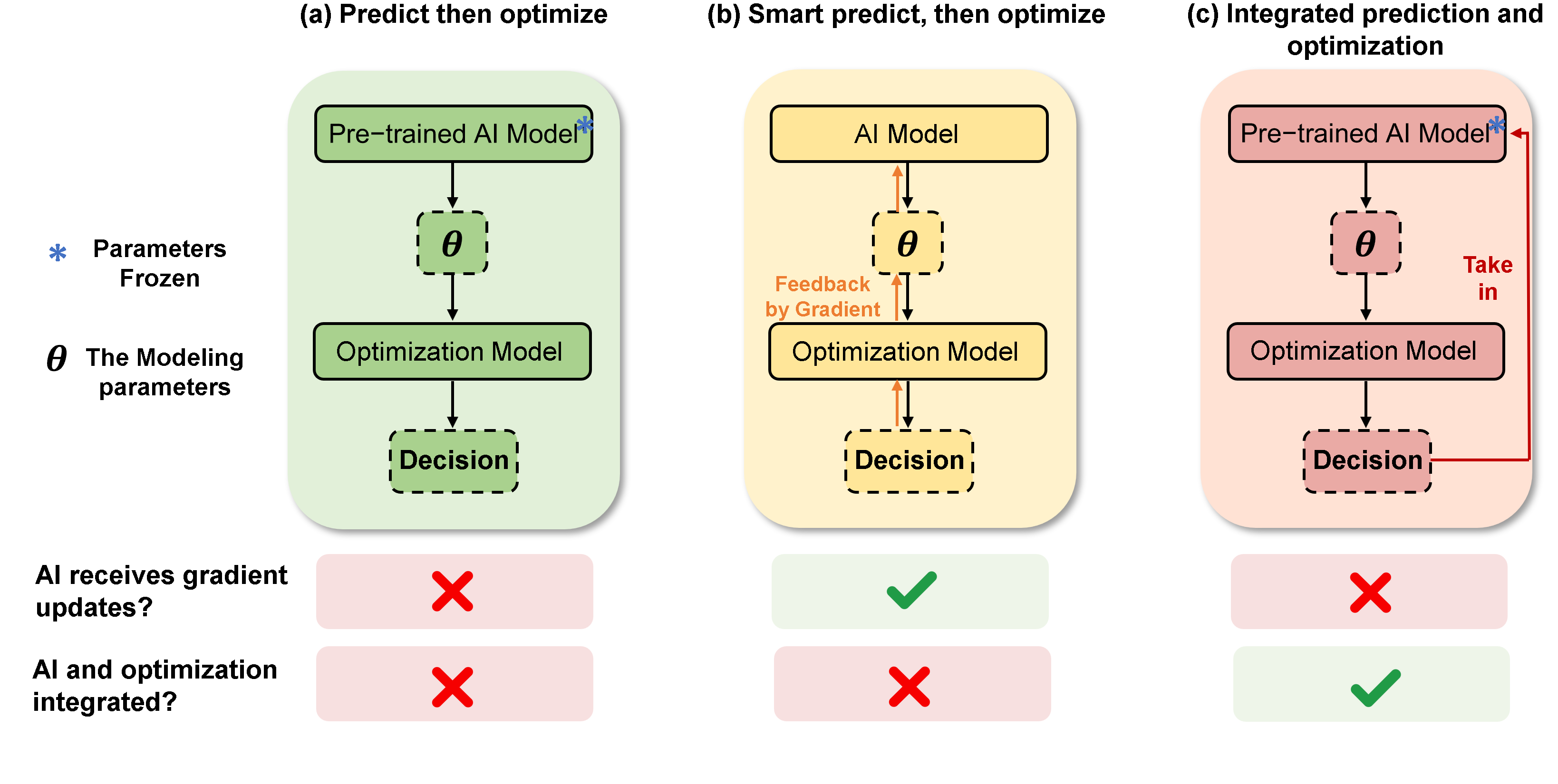}
    \caption{Illustration of how AI models generate parameters and interact with the optimization model. }
    \label{fig:param-gen}
\end{figure} 

As we discussed before, parameter generation is a crucial step of the operations research procedure. The raw data usually cannot be directly integrated into the optimization model. 
The parameter of the optimization model is traditionally obtained from human experts. 
However, it is a frequent case that some constraints and objectives cannot be described by explicit formulas. In this case, the AI model is suitable. 
Figure \ref{fig:param-gen} illustrate our three-way taxonomy of existing works according to how the AI model generates parameter and interact with the optimization model. 

\subsection{Predict, then optimize}
One straightforward approach is {predict-then-optimize}, i.e., first predict critical unknown parameters within the optimization model ({i.e., modelling parameter}) and then leverage the optimization solver to devise a decision. 
For example, consider a vehicle routing problem that requires to be solved several times on a daily basis. The raw data is the historical traveling time on the edges of a road network based on conditions like weather, time of the day, and neighbouring edges' traffic. However, as the conditions change, the current traveling time (modelling parameter) varies. 
AI-based methods can help create predictive models that estimate modelling parameters within the optimization model, by analyzing and processing large volumes of data, uncovering valuable insights that enhance the decision-making process. 
 The optimization solver then leverages these predictions to devise an optimal or near-optimal route. 
This example highlights the importance of AI in parameter generation and showcases how most solution systems tackling real-world analytics challenges benefit from the integration of both prediction and optimization.
Formally, to predict the unobserved modelling parameters $\theta$, 
we will build a training set $D$ from historical records, including $N$ pairs of attribute $x_i$ and parameter $\theta_i$, 
i.e., \[D = \{ (x_1, \theta_1), \dots, (x_N, \theta_N) \}.\] 
Here, $x_i$ is the $i$-th attribute correlated to the modelling parameter $\theta_i$ and 
 $\theta_i$ is decided in hindsight, for example, the traveling time can be decided by letting a car actually run through a road. 
Let $m(w; \cdot)$ denote an AI model for prediction, where $w$ represents the \textbf{AI model parameters}. 
In other words, we want $\hat \theta_i = m(w; x_i)$ to be a good estimation of \textbf{modelling parameter} $\theta_i$.
It should be emphasized that {AI model parameters} $w$ is different from {modelling parameter} $\theta_i$. The key difference is that $w$ is learnable while $\theta_i$ is not. 
$\theta_i$ is collected as the target value that the AI model seeks to predict given attribute $x_i$. 

The predict-then-optimize approach first trains the AI model $m$ by 
\[w^* \in \argmin_w \enspace \frac{1}{N} \sum_{i=1}^N \|m(w; x_i) - \theta_i\|^2,\]
and then solve the optimization problem 
\begin{equation} 
    v_{\hat \theta} \in \argmin_v \enspace f_{\hat \theta}(v) \enspace\text{s.t.}\enspace v \in \Cscr_{\hat \theta} 
    \enspace\text{with}\enspace \hat \theta = m(w^*; x). \label{eq:spo-opt}
\end{equation}
where $\hat\theta$ is the predictive modelling parameter, $v$ is the decision variable, $f_{\hat\theta}(v)$ is the objective to be minimized, and $\Cscr_{\hat \theta}$ is the feasible region determined by the set of constraints with predicted modelling parameter $\hat\theta$.

\subsection{Smart predict, then optimize} \label{sec:interact-gradient}
Compared to {predict-then-optimize} approach that separates prediction and optimization into two stages, 
\citet{elmachtoub2022smart} propose an end-to-end framework called {``Smart Predict, then Optimize (SPO)"} that directly trains a predictive model such that the decision error will be minimized. The motivation behind SPO is to trade off the predictive model's accuracy in the prediction stage in exchange for a near-optimal  decision. 
The error of the prediction stage is less important because {1)} each parameter in the optimization model is not equally important. {2)} these parameters may be correlated  \citep{cameron2022perils}. 

To be specific, {1)} in the prediction stage, the model inevitably makes errors due to model capacity or data quality and scale. However, take the vehicle routing problem as an example, if our main concern is the travel time on the shortest route, then we are not concerned about the lack of numerical precision in estimating the travel time for roads that are impossible to choose. 
SPO achieves better tradeoffs by allowing some error on unimportant parameters but keeping the final decision near-optimal. 
{2)} the correlated parameters are commonly seen in stochastic optimization settings. Meanwhile, it is possible that one parameter is related non-linearly to multiple prediction targets but the predictive model is unaware of the combination. For example, the predictive model predicts the cost per unit and the number of units independently but later on the total cost is used as a parameter in the optimization model. Without knowing the correlation, a small error in the prediction stage may lead to a large error in the decision. 

Formally, the general SPO framework directly aims to minimize the decision error
$\ell(v_\theta, v_{\hat \theta})$ between the two decisions from two mathematical models with ground-truth parameter $\theta$ and predictive modelling parameter $\hat\theta$.
Here, $\ell(v_\theta, v_{\hat \theta})$ could be the square of $L_2$ norm $\|v_\theta - v_{\hat \theta}\|^2$ or the objective difference $|f_\theta(v_\theta) - f_\theta(v_{\hat \theta})|$.
Recall that we only have empirical historical data $D = \{ (x_1, \theta_1), \dots, (x_N, \theta_N) \}$,
 so we can only build the empirical training loss $\frac{1}{N} \sum_{i=1}^N \ell( v_{\hat \theta_i} , v_{\theta_i})$.
This loss relies on model parameter $w$ since the predictive modelling parameter $\hat \theta_i$  relies on  $w$.
We denote this empirical loss  by $\Lscr(w) \equiv \frac{1}{N} \sum_{i=1}^N \ell( v_{\hat \theta_i} , v_{\theta_i})$ for simplicity.
The optimal model parameter $w^*$ is acquired by solving the following well-known empirical loss minimization problem 
\[w^* \in \argmin_w \enspace \Lscr(w)  \equiv \frac{1}{N} \sum_{i=1}^N \ell( v_{\hat \theta_i} , v_{\theta_i})
\enspace\text{with}\enspace \hat \theta_i = m(w; x_i)\]

The main technical bottleneck for solving this optimization problem is computing the sub-gradient 
\[\partial \Lscr(w) = \sum_{i=1}^N \left(\frac{\partial m(w; x_i)}{\partial w}\right)^\top\left(\frac{\partial v_{\hat \theta_i}}{\partial \hat \theta_i}\right)^\top
\left(
\frac{\partial \ell( v_{\hat \theta_i} , v_{\theta_i}) }{\partial v_{\hat \theta_i}}
\right),\]
which requires differentiating the $\argmin$ operator since $v_{\hat \theta_i}$ is obtained by solving an optimization problem.  
As a result, $v_{\hat \theta_i}$ can be discontinuous and non-differentiable {w.r.t} $\theta_i$. 
To overcome this challenge, the literature considers different optimization problems and accordingly propose different methods \citep{elmachtoub2022smart,amos2017optnet,surrogate-soft-constr,mandi2020smart,spo-sat,spo-black}.
Under the general SPO framework, \citet{elmachtoub2022smart} consider a constrained problem with a linear objective and a feasible region that is irrelevant to $\theta$ in Equation~$\eqref{eq:spo-opt}$. Formally, it is 
\begin{equation} 
    v_{\hat \theta} \in \argmin_v \enspace {\hat \theta}^\top v \enspace\text{s.t.}\enspace v \in \Cscr
    \enspace\text{with}\enspace \hat \theta = m(w; x). \label{eq:spo-opt-linear}
\end{equation}
Accordingly, the decision error is $\ell(v_\theta, v_{\hat \theta}) =\theta^\top ( v_\theta- v_{\hat \theta} ) $.
They propose a convex surrogate loss through the dual interpretation of $\ell(v_\theta, v_{\hat \theta})$ and obtain a subgradient for $\partial \Lscr(w)$. 
\citet{amos2017optnet} integrate the QP optimization as a differentiable layer into the AI model, by differentiating the KKT optimality conditions. 
\citet{surrogate-soft-constr} further extend it to the case that soft constraints appear in the objective. 
Soft constraints are sometimes required in practice, for example,
such constraints allow a slight excess of supply over demand through a penalty term in the objective. 
This penalty term often appears as the non-differentiable max operator function, i.e., $\max(\cdot,0)$.
This non-differentiability challenge is addressed using a proposed piecewise linear surrogate. \citet{mandi2020smart}  extend the SPO framework to NP-hard discrete optimization problems, by continuously relaxing the MILP and offering approximate subgradients. \citet{spo-black} obtain the approximate gradient by viewing the MILP as a black box and invoking two calls to an optimization solver (one with original parameter $\theta$ and another with perturbed parameter). 

It has been shown in the literature that the SPO end-to-end framework performs better than the two-stage predict-then-optimize approach. However, the two-stage approach still has some pros. 
In the two-stage approach, one may add robustness (i.e. inductive bias) in the prediction stage by observing the entire dataset or integrating common knowledge. 
This could improve the performance of the predictive model, especially with a small data size.  
It will be interesting to combine this advantage into the SPO framework in future work. 

\subsection{Integrated prediction and optimization} \label{sec:interact-loop}

In the following section, we will first describe the ``Integrated prediction and optimization'' approach through an example and then discuss the two tools developed: OptiCL \citep{mip-constraint} and  JONAS  \citep{bergman2020janos}. 

There are two motivations for this approach. 
The first motivation 
is inherited from the ``predict then optimize'', i.e.,  there may not exist an explicit formula for calculating certain parameters in the optimization model. 
The second motivation is that predicting these parameters may rely on the decision made after solving the optimization model. This additional motivation explains why, in Figure \ref{fig:param-gen}, the AI model takes in the decision variables. 
For instance, consider a telecommunication company that wants to attract customers with low prices while maximizing revenue. 
Let $x$ be the decision variable on price, $y$ be the probability of customer churn, and $\alpha$ be the additional information about customers like demographic and geographic information.  
Customer churn probability $y$ depends on both price $x$ and other information $\alpha$. 
The company collect historical dataset $D$ and trained an AI model $\hat{h}_{D}$ that can estimate $y$ given aforementioned information, i.e., $y =\hat{h}_{{D}}(x, \alpha)$. 
There may be certain policies like: 
1) the price in certain regions has a lower bound. More generally, prices $x$ fall in the feasible region $\Xscr(\alpha)$ decided by geographic information $\alpha$. 
2). If a customer is highly likely to leave then the price can be adjusted accordingly, but this adjustment varies across different regions. More generally, this variation can be expressed by an explicit formula $g$. And we write the policy as a constraint $g(x,y,\alpha) \leq 0$. 
The objective of revenue related to price $x$, customer churn probability $y$, and other information $\alpha$ is denoted by $f(x,y,\alpha)$. Then, the optimization problem is denoted as follows. 

\begin{equation}
\begin{array}{rl}
\max_{x, y} &  f(x, y, \alpha) \\
\text { s.t. } & g(x, y, \alpha) \leq 0, \\
                & y =\hat{h}_{{D}}(x, \alpha), \\
                & x \in \mathcal{X}(\alpha).
\end{array} \label{eq:integrated}
\end{equation}

Then, we will discuss how optimization problem \eqref{eq:integrated} is solved. 
The challenge is that constraint $y =\hat{h}_{{D}}(x, \alpha)$ involves a pre-trained AI model $\hat{h}_D$. 
This AI model can be linear models, decision trees, or multi-layer perceptions (MLPs). 
OptiCL \citep{mip-constraint} and JONAS \citep{bergman2020janos} are similar methods to overcome this challenge, where the constraint $y =\hat{h}_{{D}}(x, \alpha)$ are converted to linear inequalities. This conversion process is called ``embedding AI model into the optimization model''. 
1) In the case of linear models, for instance, linear regression and support vector machine, the decision region of the learned model is characterized by a half-space, thus easily converted to a linear inequality constraint. 
2) For the decision tree, since each leaf node corresponds to a polytope, the model's decision region is a set of polytopes. Conversion of the decision region to linear inequalities is a straightforward process.
3) The conversion of MLP is not as trivial. MLP often use ReLU operators, i.e., $\text{ReLU}(x)=\max(x,0)$.
This operator can be converted to linear constraints with big-M and integer artificial variables. 
Then, by recursively peeling off layers of the MLP and introducing artificial variables, the MLP with an arbitrary number of hidden layers and nodes can be embedded into the optimization model. 

While OptiCL and JANOS are similar regarding the conversion process, they differ in some of the detailed techniques. 
JANOS employs a discretization method that breaks down the sigmoid operator, $\sigma(x) = \frac{1}{1 + e^{-x}}$, into piecewise linear functions.  
OptiCL proposes two techniques for handling the uncertainty of the trained AI model. The uncertainty arises because the true constraints are unknown and the AI model may not capture the true constraints accurately. 
To mitigate these challenges, OptiCL adopts the following strategies: 1) Employs an ensemble of AI models. A solution is deemed feasible only if a significant portion of the ensemble agrees on it.
2) Defines the dataset's convex hull as a trust region. Predictions are deemed reliable when they are within or around this trust region. Additionally, 
recall that the modeling parameter $\hat\theta$ is predicted by AI model $m(w; \cdot)$ from observed attributes $x$, denoted by $\hat \theta = m(w; x)$. This process is commonly termed as point estimation. It fails to consider the uncertainty of $\hat \theta$. The uncertainty exists because we cannot observe every related attribute and the AI model cannot perfectly predict the target modeling parameter.  \citet{bertsimas2020from} propose accounting for the uncertainty when generating the model parameters. They also extend their approach to incorporate the interaction between prediction and optimization. The authors apply the proposed approach to an inventory management problem.

\citet{opt-constraint} provide a survey on ``optimization with constraint learning''. Though with a different name, it actually presents predictive AI models that are integrated into constraint and objective.  
\citet{boosting} present another survey on the same topic but focus on discrete optimization problems. 
They cover a broad range of methods where AI models are integrated into optimization models, including: 
1) Building explicit constraint equations using active learning. Active learning is an iterative process involving a dynamic interaction between domain experts and a learning model. Initially, the learning model formulates constraints and employs a solver to derive a solution. The domain expert then assesses the feasibility of this solution. Feedback from this assessment refines and enhances the learning model's understanding and subsequent constraint formulations. 
2) Approximate the objective using simple formulas. In scenarios where evaluating the objective function is computationally expensive, such as involving simulations, AI models are employed to approximate the objective. The aim is to capture the internal functions of the simulator through this approximation, thereby reducing computational overhead but still finding a high-quality approximate solution.



To conclude, the AI model is able to generate key parameters for the optimization model, eliminating manual work and surpassing human capabilities, specifically where explicit formulas are not applicable. 
We discussed three categories of how the AI model interacts with the optimization model. A  direction left unexplored in the literature involves allowing AI models to receive gradient feedback {while} integrating AI with optimization. 
The motivation for this direction is that at the beginning, we do not know how to formulate a certain constraint, so we randomly initialize it as an AI model. Obviously, the decision made by solving the optimization model will be random. However, the infeasibility or error of decision serves as feedback to improve the AI model. 
The challenge of realizing this direction lies in the computation cost. 
Embedding even the simplest MLP can introduce many additional (integer) artificial variables, consequently increasing the scale of the optimization problem. 
If the scale of the optimization problem goes beyond a manageable limit, the gradient of the decision error w.r.t AI model parameter will not be calculated within an acceptable time. 
Nevertheless, it may remain promising to explore this direction starting from a small-scale application. 

\section{Model Formulation} \label{sec:model_formulation}

In the operations research process, modelling holds paramount importance and is often a time-consuming step. Modelling is a defining characteristic of OR, and thus, it warrants significant attention. Modelling in OR often involves the creation of mathematical models, which typically consist of three main elements: decision variables, constraints, and objective function(s). Decision variables are employed to represent specific actions under the control of the decision-maker. 
In complex models, it is common to define auxiliary or artificial variables to better model the problem. Although not directly controlled by the decision-maker, these variables are also considered decision variables. Constraints serve to establish limits on the range of values that each decision variable can assume, with each constraint typically translating a specific restriction (e.g., resource availability) or requirement (e.g., meeting contracted demand) within the model. Constraints dictate the feasible values assignable to decision variables, effectively determining the possible decisions for the system or process under consideration. The final component of a mathematical model is the objective function, which represents a mathematical expression of a performance measure (e.g., cost, profit, time, revenue, or utilization) as a function of the decision variables. Regarding the nature of the objective function, the goal is typically either to maximize or minimize its value.

In recent years, the most significant breakthrough in AI research has been the remarkable progress in Natural Language Processing (NLP) achieved by large language models (LLMs). There have been preliminary efforts to explore the potential of using LLMs, such as ChatGPT \citep{brown2020language} and Llama \citep{touvron2023llama}, for mathematical modeling. These early investigations have demonstrated promising results, showcasing the potential for LLMs to revolutionize the process of formulating mathematical models in operations research by automating the translation from natural language problem descriptions into mathematical representations.

Large language models are obtained through a complex, multistage process \citep{ouyang2022training,zhou2023lima,alpaca}. It starts with pre-training, where the model is exposed to a vast amount of text data from a broad spectrum of domains. This enables the model to learn general language patterns, semantics, syntax, and common-sense knowledge.
During pre-training, the model attempts to predict the next word in a sentence, thus learning the context-dependent representations of words. Once this unsupervised learning stage is complete, the model undergoes fine-tuning. In fine-tuning, the model is instructed with a more specific task using a smaller, labeled dataset. This step helps the model refine its capabilities and align its responses with specific goals. Some LLMs \citep{ouyang2022training} also incorporate a level of reinforcement learning from human feedback (RLHF), where the model's responses are ranked by human preference and used to further improve its predictions.

When utilized as a backend model during inference time, LLMs offer impressive capabilities. One important application is reference editing \citep{gao-etal-2023-reference,reid2022learning,wen2017network}, where the LLM is used to generate, evaluate, and modify drafts of text based on user input. In this process, the model first generates an initial draft. Then, based on feedback or additional instructions, the model performs iterative edits until the text meets the desired quality and specifications. This not only enables the creation of more refined and context-specific content but also provides a mechanism for users to control the model's output in a more flexible manner. Through such interactive usage, LLMs can be a powerful tool in various fields, from content creation and editing to customer service and beyond.

In this survey, we aim to harness the natural language understanding capabilities of large language models to generate mathematical models based on requirement lists or problem descriptions provided by business owners. Additionally, these models can be modified by experts through natural language inputs, streamlining the model iteration process and fostering a more collaborative approach. 
In this section, the Llama family LLMs \citep{touvron2023llama, rozière2023code}, considered as one of the most powerful open-sourced LLMs, is discussed and adopted for generating mathematical models from natural language problem descriptions.
We analyze problems with two different levels of complexity as described in sections \ref{sec:textbook_modelling} and \ref{sec:real_modelling} and we provide insights into the effectiveness of using LLMs for the model formulation stage.


\subsection{Textbook modelling problems} \label{sec:textbook_modelling}


\begin{figure}[t]
    \centering
    \includegraphics[width=\textwidth]{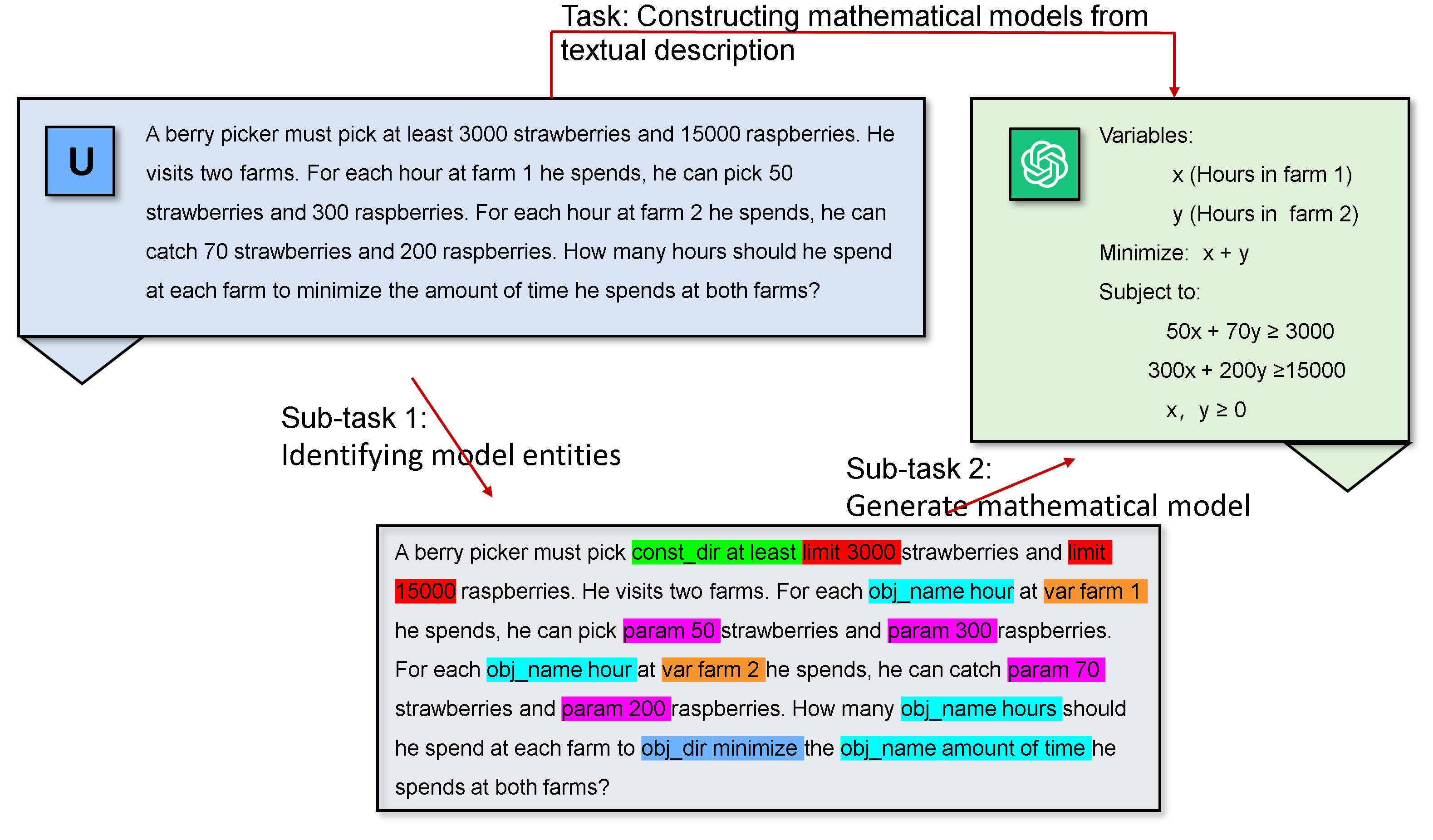}
    \caption{Example on a textbook modeling problem.}
    \label{fig:model_ex1}
\end{figure} 

In this section, we focus on the textbook modeling problem. The problems of this level are quite simple. Figure \ref{fig:model_ex1} demonstrates an example. 
We utilize the dataset from the NL4OPT competition held at NeurIPS 2022 \citep{ramamonjison2023nl4opt}, which was centered around constructing mathematical models from natural language problem descriptions.  This dataset contains 713 training, 99 validation, and 289 testing data points.  In this dataset, each data point consists of both a problem description and a human-composed mathematical model. Declaration-level mapping accuracy is proposed in this competition for evaluating various methods. As shown in Figure \ref{fig:model_ex1}, the competition decomposes the overall task into two sub-tasks: identifying model entities (sub-task 1) and generating the mathematical model (sub-task 2). For an in-depth understanding of the competition and its evaluation criteria, we refer readers to \citep{ramamonjison2023nl4opt}. 
{The purpose of using the NL4OPT dataset here is to evaluate the performance of LLM for text-book level problems, review popular LLMs, and provide some insights.}

As mentioned before, most LLMs are obtained through a two-step process: pretraining and finetuning. During pretraining, LLMs are exposed to large volumes of text data, enabling them to grasp language patterns and common knowledge. The finetuning stage then refines the model on specific, high-quality instruction-answer pairs, ensuring it adheres to given instructions. 
Training over-parameterized models on a vast amount of text data leads to the emergence of advanced capabilities in LLMs, such as contextual comprehension, generating coherent responses, and even solving mathematical problems. These skills, not explicitly taught, naturally arise from the model. This phenomenon is known as ``emergence".

The usage of LLM has become increasingly standardized. Developers often download a pre-trained LLM of their chosen size and finetune it with specific data to meet their needs. Various off-the-shelf LLMs are obtained in this way. Their difference lies in model size (the amount of disk space the parameters of the model take) and finetuning dataset. To showcase the superiority of LLMs in modeling textbook problems and to give relevant insights, we investigated several widely-used off-the-shelf LLMs, including Llama-2-13b-chat \citep{touvron2023llama}, Code-Llama-34b-instruct \citep{rozière2023code} and Llama-2-70b-chat \citep{touvron2023llama}, a variant LLM further fine-tuned with the NL4OPT training set, named ``Llama-2-13b-chat (SFT)”, and the NL4OPT winning submission for sub-task 2. 
Here, the winning submission also follows a similar structure, i.e., pick an off-the-shelf pre-trained LLM -- Bart \citep{bart} and fine-tune it with the NL4OPT train set. A careful hyper-parameter tuning is conducted on the validation set to achieve high performance. 
The same metric (Declaration-level mapping accuracy) is employed to evaluate different LLMs (See \cite{ramamonjison2023nl4opt} for more details). 
Table \ref{tab:nl4opt-res} lists the performance of different LLMs and the NL4OPT winning submission. 
\begin{table}[tbp]
\centering
\scalebox{.97}{

\begin{tabular}{c|l | r|l | c}
\toprule
    Index & \textbf{LLMs} & Model size & Finetuning dataset & \textbf{Acc.}  \\ \midrule
   1& Llama-2-13b-chat  & 52 GB & Open domain & 24\% \\ 
    2&Code-Llama-34b-instruct  & 136 GB & Programming and math domain & 65\% \\ 
   3& Llama-2-70b-chat  & 280 GB & Open domain & 37\% \\ 
   4& Llama-2-13b-chat (SFT)  & 52 GB & Open domain+NL4OPT train set  & 82\% \\ 
  5&  NL4OPT winning submission & 1 GB& NL4OPT train set & 90\% \\ \bottomrule 
\end{tabular}

}
\caption{\label{tab:nl4opt-res} 
Declaration-level mapping accuracy (denoted by ``Acc.") for different LLMs and NL4OPT Winning Submission. 
}
\end{table}

Without any training, Code-Llama-34b-instruct has achieved 65\% test accuracy.
Without hyperparameter tuning, Llama-2-13b-chat (SFT) achieves an 82\% accuracy. 
This already concludes that LLM is easy to use and performs well on textbook-level problems. 
The winning submission archives the 90\% performance, but we need to mention that it is unfair to compare the first four LLMs in Table~\ref{tab:nl4opt-res} with the winning submission. 
The first four LLMs solve the task of modeling from natural language in an end-to-end manner, 
while as mentioned before, the winning submission of the competition decomposes the overall task into two sub-tasks, thereby complicating the process. 
The reported 90\% performance for sub-task 2 assumes that the entity information is perfectly annotated after sub-task 1, indicating that this 90\% is merely an upper bound and does not measure the end-to-end performance. 
Lastly, it is worth noting that LLMs can be easily further improved. For example, with a more comprehensive and high-quality pretraining corpus, such as operational research papers, LLM will understand and generalize better for mathematical modeling. 
When comparing the performance of Code-Llama-34b-instruct and Llama-2-70b-chat, it becomes evident that a larger model size does not necessarily guarantee improved performance. The quality and relevance of the fine-tuning dataset play a more significant role. Further, by examining Llama-2-13b-chat against its further fine-tuned version, we reinforce the importance of relevant and in-distribution data for a fine-grained task like modeling from natural language. 

To summarize, the current performance of LLMs for modeling textbook-level problems has been impressive. 
As these LLMs continue to evolve and improve, it is expected that their capability for effective modeling directly from natural language inputs will significantly advance in the near future.

\subsection{Real-world problems} \label{sec:real_modelling}

\begin{figure}[t]
    \centering
    \includegraphics[width=\textwidth]{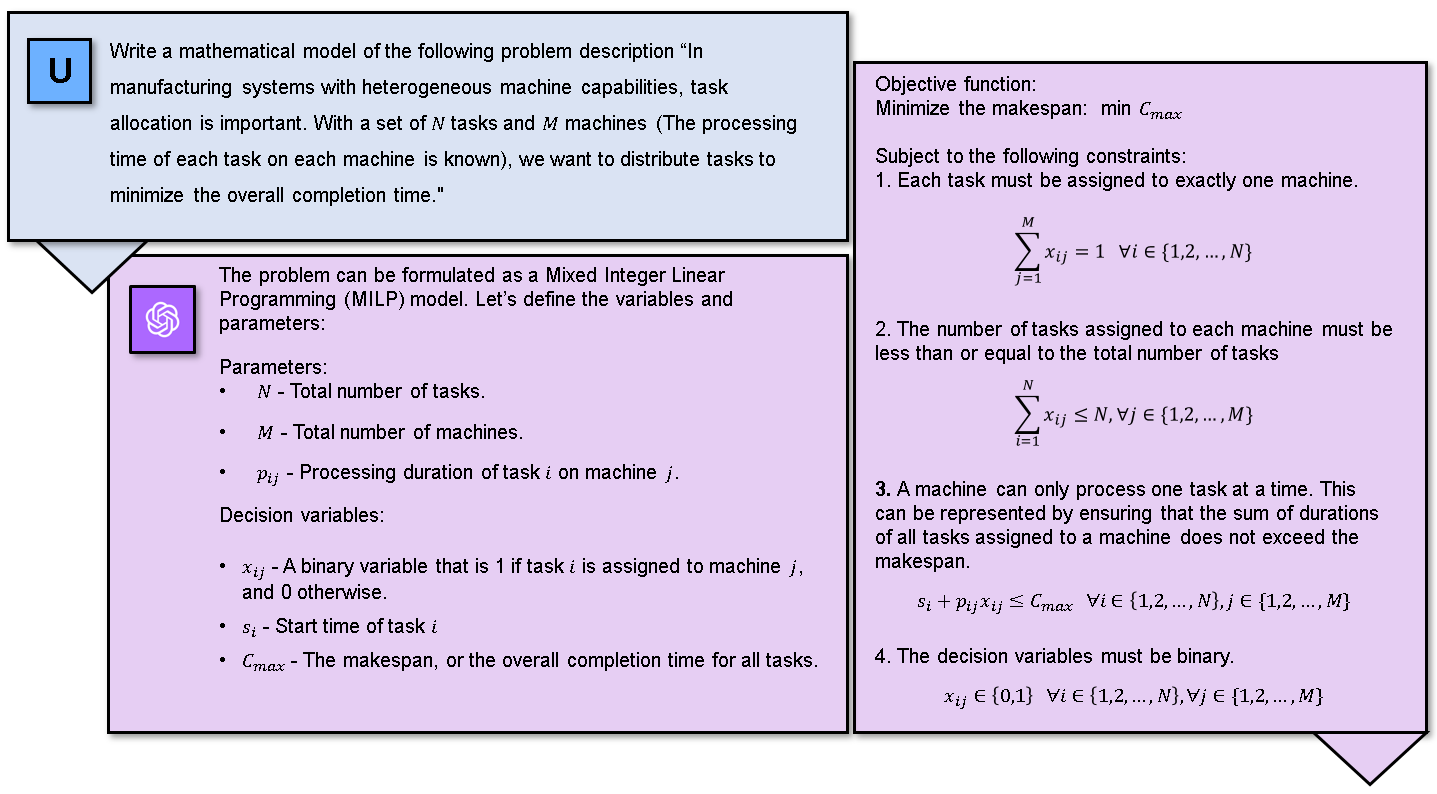}
    \caption{Example on the Unrelated-machines Scheduling Problem.}
    \label{fig:model_ex2}
\end{figure}

In the second experiment, we use a set of real-world problems, including Unrelated-Machine Scheduling Problem \citep{BLAZEWICZ1991283}, One-dimensional Bin Packing \citep{pinedo2016scheduling}, Capacitated Vehicle Routing \citep{toth2014vehicle}, Team Formation \citep{anagnostopoulos2012online}, Portfolio Optimization \citep{markowitz1952portfolio}, Staff Scheduling \citep{blochliger2004modeling}, and Airline Revenue Management \citep{talluri2004revenue}, Cutting Stock Problem in the Printing Industry~\citep{cutting}, Fair Distribution of Relief Aid Supplies in Disaster Response~\citep{relief}, and Short-term Scheduling in Mining Operations~\citep{mine_short}. 
For each problem, we provide a detailed problem description and use the simple prompt ``Write a mathematical model of the following problem description + \textless{}Problem Description\textgreater{}" as input to the Llama-2-70b-chat LLM. We show the example of the Uniform Machine Scheduling Problem \citep{BLAZEWICZ1991283} in Figure~\ref{fig:model_ex2}. 
The generated model formulation, although it contains several flaws, offers a starting point for further refinement.
On the constructive side, LLM correctly identifies sets, parameters, and important decision variables, and constructs a skeletal framework for further refinement. 

However, the generated formulation still involves errors that need to be addressed, e.g., incorrect constraint expressions, extra constraints, and missing constraints. These errors can render the problem infeasible or the final solution irrelevant to the problem.
For the specific example in Figure~\ref{fig:model_ex2}, Constraint 2 is actually redundant due to the existence of Constraint 4. 
Furthermore, a decision variable $s_i$ is introduced, denoting the start time of task $i$, without accompanying constraints to prevent the temporal overlap of tasks assigned to the same machine. 
However, if we point out the redundant or incorrect constraints and prompt Llama to correct the errors, it can make the necessary corrections after a few rounds of prompts and responses. 
For the above-mentioned example, 
LLM finds that introducing $s_i$ is unnecessary in the context of the Unrelated-Machine Scheduling Problem because the order of tasks on the same machine does not influence the machine's completion time. 
After removing $s_i$, it suggests an improved Constraint 3: $\sum_{i=1}^{N} p_{ij} x_{ij} \leq C_{\text{max}}, \quad \forall j \in \{1,2,\ldots,M\}$.
This demonstrates the potential for iteratively refining the generated mathematical models through interaction with LLM, leading to more accurate and relevant formulations.

Overall, while LLM made several errors, the provided formulations can serve as a starting point for OR experts to create mathematical models. However, OR experts should not rely on LLM to accurately create mathematical models, especially for less common or complex problems. Each output needs to be thoroughly verified and adjusted by the experts to ensure correctness and relevance.

\section{Automatic Algorithm Configuration} \label{sec:solver_preparation}

\begin{figure}[ht]
    \centering
    \includegraphics[width=0.67\textwidth]{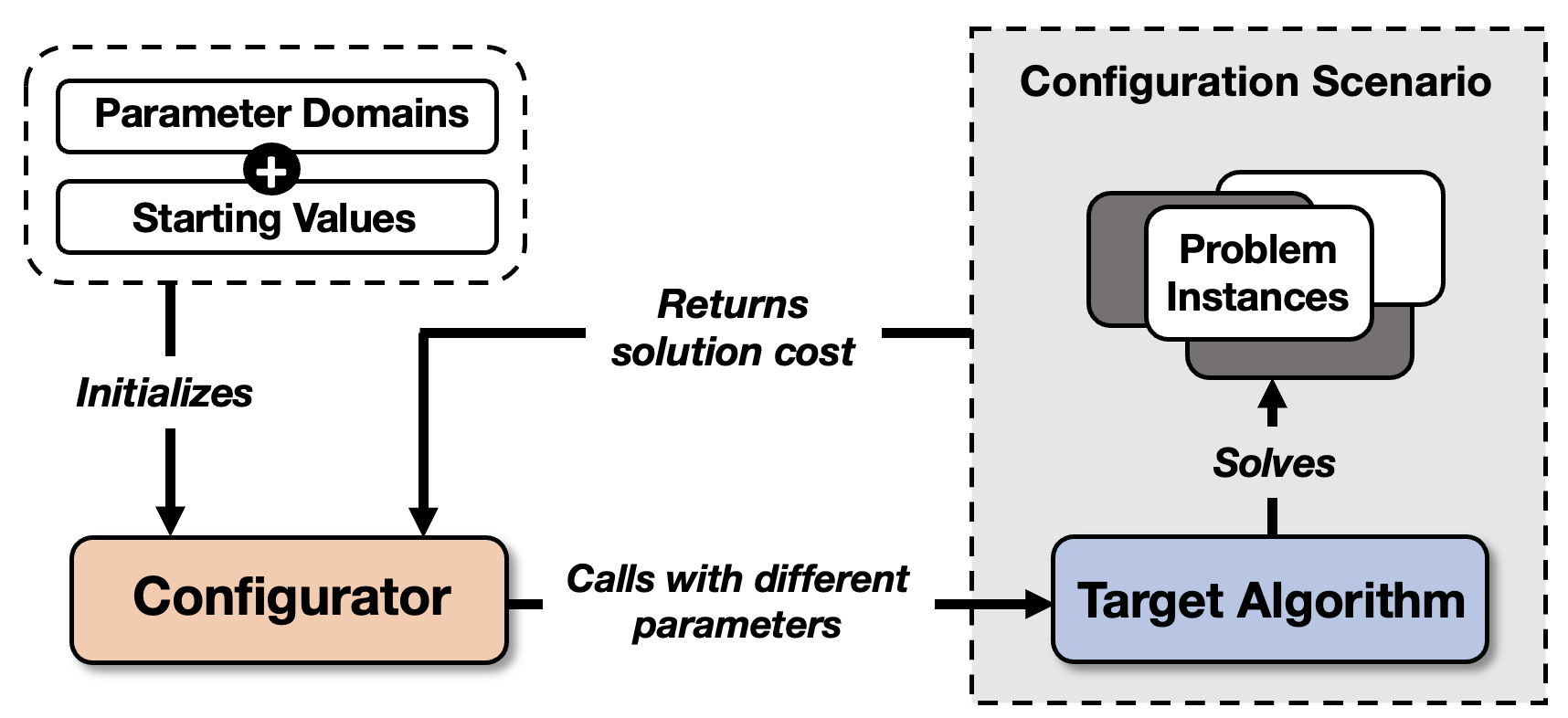}
    \caption{Illustration of automatic algorithm configuration, adapted from the work by  \citet{hutter2009paramils}. }
    \label{fig:auto-conf}
\end{figure}

Optimization solvers contain parameters that significantly impact their performance. However, manually tuning these parameters is challenging due to their complex interactions. Automatic algorithm configuration (AAC) aims to automate the tuning process by systematically exploring the parameter space to find optimal configurations that maximize performance. More formally, consider a parameterized optimization algorithm or solver $\Ascr$, its parameter space $\Theta$, a distribution of problem instances $\Dscr$, and a cost function $o$ that evaluates the performance of a solver configuration $\theta \in \Theta$ on an instance $d$. The goal of AAC is to solve:
\begin{align} \label{eq:AAC}
    \theta^* &= \argmin_{\theta \in \Theta} \mathop{\mE}_{d \sim \Dscr} o(\Ascr(\theta), d)
\end{align}
to find the configuration that optimizes expected performance over $\Dscr$.

The Parametric Iterated Local Search (ParamILS) framework \citep{hutter2009paramils} is a representative AAC method for tuning optimization solvers. 
Given a training set of instances $\{d_1, ..., d_n\}$, ParamILS proceeds as follows:
\begin{itemize}
    \item {\bf Initialization:} Randomly sample $r$ configurations $\{\theta_1, ..., \theta_r\}$ and select $\theta$ with the lowest training cost:
    \[\theta = \argmin_{\theta \in \{\theta_1, ..., \theta_r\}} \Lscr(\theta) =  \frac{1}{n}\sum_{i=1}^n o(\Ascr(\theta), d_i)\]
    \item {\bf Iterated  Local search:} While the stopping criteria (e.g., solution quality exceeds a threshold, maximum runtime is reached) is not met, {local search} and {perturbation} are repeatedly executed. 
    \begin{itemize}
        \item {\bf Local search: } Randomly sample a new configuration $\theta'$ from the neighborhood $ \Nscr(\theta)$. If $\Lscr(\theta') < \Lscr(\theta)$, set $\theta = \theta'$. Repeat this for $k$ times. Here $k$ is a hyperparameter.
        \item {\bf Perturbation:} Randomly choose a configuration $\theta'$ differing from $\theta$ in at most $p$ parameters where $p$ is a pre-defined hyperparameter. Perform the local search from $\theta'$. If better, set $\theta = \theta'$. 
    \end{itemize}
\end{itemize}

ParamILS has shown significant benefits for configuring mixed integer programming solvers and other optimization algorithms. More specifically, in \cite{hutter2010automated}, the authors apply ParamILS to leading MIP solvers: CPLEX~\citep{cplex2009v12}, Gurobi~\citep{gurobi}, and LPSOLVE \citep{berkelaar2015package}. 
To assess the benefits of auto-tuning, the authors perform experiments on two sets of MIP instances: a library of over 200 cases from prior studies, and over 90 new larger cases. They measure the runtime required for the solvers to prove optimality or find the best solution within a time limit. Results show that auto-tuning reduces the runtime of CPLEX and Gurobi by over 25-50\% on average compared to the default configurations, and over 50-100\% for the hardest problem instances. The authors also analyze how the parameter settings found by ParamILS differ based on instance features like the number of variables and constraints. They find that more aggressive parameter configurations (e.g. stronger presolving) are selected for easier instances, while more robust and diverse settings are chosen for harder cases. This demonstrates how auto-tuning can exploit instance properties to configure solvers more appropriately.

Alternatively, \citet{lopez2016irace} propose an R package called IRACE for automatic algorithm configuration. IRACE implements a method called iterated racing for efficiently exploring the parameter space. In one iteration of IRACE, multiple parameter configurations are evaluated in parallel on training instances with limited budgets (called racing), and poor performers are discarded. 
The purpose of racing is to gather preliminary performance and estimate the configurations' quality. 
This process is repeated over multiple iterations and strikes a balance between exploring and exploiting the promising regions of the parameter space. Compared with ParamILS, which handles mainly categorical parameters, IRACE handles a wide range of parameter types, including continuous, integer, categorical, and conditional parameters. Meanwhile, IRACE leverages statistical testing to select candidate configurations, and thus is more suitable when the performance of the algorithm being tuned exhibits stochastic behavior. 
\citet{lopez2016irace} evaluate IRACE on configuring algorithms from various domains, including constraint programming and hyperparameter optimization. Results show that IRACE can achieve average speedups of 20-60\% over default parameter settings and matches or outperforms manual tuning by experts. For algorithms with a large number of parameters, the benefits of auto-tuning with IRACE become even more substantial. 

In Equation~\ref{eq:AAC}, the expected performance on a set of instances is called the performance function, denoted as 
$g(\theta) = \mathop{\mE}_{d \sim \Dscr} o(\Ascr(\theta), d)$. 
However, when navigating through high-dimensional parameter space, the curse of dimensionality makes both random samplings \citep{lopez2016irace} and local search strategies \citep{hutter2009paramils,hutter2010automated} often inadequate.
To capture the knowledge of parameter space, building a surrogate performance function is advantageous. 
Model-based methods \citep{seq-model-aac,Lindauer_Hutter_2018,JMLR:v23:21-0888} incorporate a learned model as the surrogate performance function to guide the search for configurations. It alternates between exploiting the current model to find promising configurations and exploring the configuration space to improve the model. However, it fails when dealing with stochastic or non-smooth functions $g$, e.g., multiple evaluations at the same configuration can result in different outcomes for randomized algorithms. 
More future research is needed to improve the efficiency of exploring the configuration space, e.g., 
using more advanced sampling methods, combining the idea of racing, sampling, and surrogate model \citep{f-race-sampling,AnaHoo20b,AnaHoo20}.

To reiterate, AAC automatically finds the best-performing configuration settings for any given algorithm, with respect to a particular set of problem instances. 
In essence, AAC targets a {single} algorithm and treats it as a ``{black box}". 
This carries both merits and drawbacks. On the upside, the algorithm is treated as a black box that accepts an algorithm configuration and an optimization problem, then yields performance metrics, such as solving time. 
It means there's no need to comprehend the internal workings of the algorithm and empowers {non-experts} to enhance the performance of optimization algorithms. Moreover, by exploring the parameter space in a structured yet randomized manner, AAC might outperform experts by discovering complex parameter interactions that would be difficult to identify manually. However, there are two inherent disadvantages. 
Firstly,  the ``No Free Lunch" theorem \citep{wolpert1997nfl} suggests that certain algorithms may outperform others on a specific set of problem instances.
Thus, relying solely on {one} algorithm may be insufficient if the distribution of problem instance $\Dscr$ is too diverse.
Secondly, AAC's approach means it doesn't capitalize on insights derived from the algorithm's internal mechanics. 
To overcome the limitations, some closely related tasks on hyperparameter tuning are proposed. 

To address the first limitation, it is important to choose an optimization algorithm that aligns well with the nature and distinct attributes of the problems at hand. 
The task of algorithm selection is highlighted by \citet{bischl2016aslib}. 
This task involves determining the most suitable algorithm among a set of algorithms for a given problem instance
in order to exploit the varying performance of algorithms over a diverse set of instances. 
Formally, 
given $K$ candidate algorithms $\{\Ascr_1, ..., \Ascr_K\}$, a distribution of problem instances $\Dscr$, and a cost function $o$ that evaluates the performance of an algorithm on an instance $d$. 
The objective is to identify a mapping $f(\cdot)$. This mapping assigns each problem instance to an algorithm index, optimizing the following expected performance:
\begin{align} 
    f^* &= \argmin_{f} \mathop{\mE}_{d \sim \Dscr} o(\Ascr_{f(d)},d).
\end{align}
Here, 
$f(d)$ retrieves the index of the algorithm best suited for the problem instance $d$.
While the algorithms are still treated as black boxes, the characteristic inside the instances is utilized, since the mapping $f$ is often implemented by using instance features like the density of the constraint matrix. These instance features are then mapped to an algorithm via an AI model. 
It is worth noting that algorithm selection disregards algorithm configuration. Therefore, if one algorithm underperforms compared to another, it might be due to suboptimal configuration rather than the algorithm's inherent unsuitability for the task. In future research, it would be intriguing to investigate tasks that integrate both algorithm selection and automatic configuration tuning. 

Regarding the second limitation, i.e., AAC perceives the algorithm as a ``black box", one potential strategy is to group algorithms based on their inherent characteristics, thus peeling back a layer of this ``black box". 
For example, many algorithms are iterative and can benefit from dynamic configuration during their execution, adjusting based on the information available at run time.
\citet{adriaensen2022automated} propose a general automated {dynamic} algorithm configuration framework for iterative algorithms. 
Another promising approach is to focus on a single specific algorithm, fully ``unboxing" it to leverage its inherent features. Advanced AI-driven iterative optimization algorithms are further explored in Sections~\ref{sec:learn4continuous} and \ref{sec:branch_bound}. 

\section{Algorithm Selection and Design for Continuous Optimization} \label{sec:learn4continuous}

In this section, we review the literature on using AI techniques to enhance the algorithms for solving continuous optimization problems. 

\subsection{Enhancement for gradient-based methods} \label{sec:first_order}
In this section, we consider solving an unconstrained continuous optimization problem
\begin{equation} \label{prob:gradient}
    \min_{x \in\mR^n}\enspace f(x),
\end{equation}
where $f:\mR^n\to\mR$ is a differentiable function. Gradient descent~\citep{beck2017first} is one of the most widely used algorithms for solving Problem~\eqref{prob:gradient} because of its cheap cost at each iteration. However, the performance of gradient descent is quite limited by the fact that it only makes use of the latest gradient and ignores past information. To resolve this issue, many gradient-based optimization algorithms have been proposed to improve the performance of gradient descent, and we summarize a few representatives below:
\begin{itemize}
    \item Gradient descent 
    \begin{equation*}
        x^{t+1} = x^t - \eta \nabla f(x^t),
    \end{equation*}
    where $\eta$ represents the learning rate.
    
    \item Gradient descent with momentum~\citep{tseng1998incremental}
    \begin{equation*}
        \begin{split}
            g^{t} &= \gamma g^{t-1} + (1-\gamma)\nabla f(x^t), \\
            x^{t+1} &= x^t - \eta g^{t},
        \end{split}
    \end{equation*}
    where $\gamma \in [0, 1)$ is the momentum coefficient determining the effect of the previous gradient updates on the current update. 
    \item Nesterov's accelerated gradient descent~\citep{Nesterov1983AMF}
    \begin{equation*}
        \begin{split}
            g^{t+1} &= \gamma g^t - \eta \nabla f(x^t + \gamma g^t), \\
            x^{t+1} &= x^t + g^{t+1}.
        \end{split}
    \end{equation*}
    \item Adaptive gradient descent (AdaGrad)~\citep{duchi2011adaptive}
    \begin{equation*}
        \begin{split}
            G^t &= G^{t-1} + \mathop{\mathbf{diag}}(\nabla f(x^t))^2, \\
            x^{t+1} &= x^{t} - \eta * \left[ (G^t)^{-\frac{1}{2}} \nabla f(x^t) \right],
        \end{split}
    \end{equation*}
    where $G_t$ is a rough approximation of the diagonal elements of the Hessian matrix and the symbol $*$ is multiplication between a scalar and a vector.  
    \item Adam~\citep{kingma2014adam}
    \begin{equation*}
        \begin{split}
            g^t &= \gamma_1 g^{t-1} + (1-\gamma_1)\nabla f(x^t),\\
            G^t &= \gamma_2 G^{t-1} + (1-\gamma_2)\mathop{\mathbf{diag}}(\nabla f(x^t))^2, \\
            x^{t+1} &= x^{t} - \eta * \left[ (G^t)^{-\frac{1}{2}} g^t \right],
        \end{split}
    \end{equation*}
    where $\gamma_1, \gamma_2 \in [0, 1)$ are momentum coefficient and squared gradient coefficient respectively. More specifically, the squared gradient coefficient $\gamma_2$ determines the importance of previously squared gradients when performing exponential moving average for $G^t$. 
\end{itemize}

We refer interested readers to \citet{ruder2016overview} for a more detailed overview of those gradient-based methods with explicit update rules. In this section, we want to review the papers aiming to learn a parameterized update rules.

As the optimization process can be regarded as a trajectory of iterative updates, LSTMs 
(see Section~\ref{sec:rnn}) are natural modelling choices for learning the update rule. \citet{andrychowicz2016learning} proposed the first LSTM-based learning technique to replace the explicit update rules in those gradient-based methods. The main idea is to learn a good update rule which is expected to have a good performance on some target problem sets. More specifically, it aims to solve the following problem 
\[
    \theta^* \in \argmin_{\theta} \mathop{\mE}_{(f, x^0)\sim \Fscr}\left[ \sum_{t=1}^T f(x^t) \right]
    \enspace\text{with}\enspace
    x^{t+1} = x^t - g^t
    \enspace\text{and}\enspace
    \begin{bmatrix}
        g^t \\
        h^{t+1}
    \end{bmatrix} 
    = m(\theta; \nabla f(x^t), h^t),
\]
where $T$ is some pre-determined maximal number of iterations, $\Fscr$ is the target problem set containing problem instances and corresponding initial points, $m(\theta; \cdot, \cdot)$ is the LSTM model with $\theta$ being the learnable parameters, and $h^t$ is the hidden embedding for all the gradient information up to iteration $t$. There are two technical difficulties with this approach. The first one is the varying dimensions, as the problem instances $f$ in the target problem set $\Fscr$ may have different dimensions, which requires the model $m(\theta; \cdot, \cdot)$ to be able to handle varying input dimensions.  \cite{andrychowicz2016learning} resolved this technical difficulty by letting the LSTM model operate coordinatewise on variable $x$, i.e., 
\[
\begin{bmatrix}
        g^t_i \\
        h^{t+1}_i
    \end{bmatrix} 
    = m(\theta; \nabla f(x^t)_i, h^t_i), \enspace \forall i = 1, \dots, n.
\]
In this way, the LSTM model can handle problem instances with any dimension. 

The second challenge is the choice of the maximal number of iterations $T$. A small $T$ may yield unsatisfactory results, which is also known as truncation bias, and a big $T$ may cause gradient explosion~\citep{pascanu2013difficulty,chen2021learning}. Much research has been done to solve this problem. \citet{lv2017learning} proposed to add random scaling and convex regularizers to stabilize the training process so that a larger $T$ can be selected. They show that their strategies are intended to prevent significant random updates when the LSTM optimizer is insufficiently trained. Alternatively, \citet{chen2020training} utilized techniques from curriculum learning~\citep{bengio2009curriculum} to gradually increase $T$ during the training so that the model can mitigate the dilemma between truncation bias and gradient explosion. 
\citet{chen2020learning} tackled this problem by learning another variational stopping policy so that $T$ is no longer a fixed pre-determined parameter. Such a design will cause difficulty in training the model. Therefore, the authors propose a novel training procedure that decomposes the task into an oracle model learning stage and an imitation stage. \citet{metz2019understanding} suggested using large $T$ and tried to overcome the gradient explosion problem by computing the gradient as the weighted average of two unbiased gradient estimators. The two unbiased gradient estimators are computed by the reparameterization trick~\citep{kingma2013auto} and the log-derivative trick~\citep{williams1992simple}. Alternatively, \citet{wichrowska2017learned} suggested using more-advanced RNN structures. In contrast to a single RNN layer, the authors implemented three RNN layers in a hierarchical manner, denoted as ``bottom RNN", ``middle RNN", and ``upper RNN". The ``bottom RNN" takes the scaled gradients as input and outputs the hidden states, the ``middle RNN" takes these hidden states and produces averaged hidden states, and the ``upper RNN" receives the averaged hidden states. The authors showed empirically that such a hierarchical design leads to lower memory and computing overhead while achieving superior generalization. 

Besides LSTM, reinforcement learning is another widely adopted framework for learning gradient-based methods. \citet{li2016learning} proposed the first RL-based framework for learning gradient-based methods. In their setting, the state $s^t$ consists of the current iterate $x^t$ and feastures $\Phi^t$ which depends on the history of iterates $x^0, \dots, x^t$, gradients $\nabla f(x^0), \dots, \nabla f(x^t)$ and objectives $f(x^0), \dots, f(x^t)$. The action $a^t$ is the step $\Delta x$ that will be used to update the iterate. The reward is defined as the decrease in the objective. Finally, they use RL to learn a policy $\pi$ that can sample the update steps from a Gaussian distribution whose mean and variance are learnable parameters. Based on this framework, \citet{li2017learning} developed an extension that is suited to learning optimization algorithms for high-dimensional stochastic problems. More specifically, they notice that if the values of two coordinates in all current and past gradients and iterates are identical, then the step vector produced by the algorithm should have identical values in these two coordinates. Based on this finding, they proposed grouping coordinates under permutation invariance into a coordinate group. 
They show that this formulation reduced the computation cost of training neural networks on well-established image classification datasets.

Table~\ref{tab:gradient} summarizes the abovementioned methods into two categories. 
There is no clear winner between LSTM-based and RL-based methods. 
When training loss is differentiable w.r.t. the model parameter, supervised learning and a medium-size model (like LSTM) are favourable, due to data efficiency and easier hyper-parameter tuning. This is the reason why more literature follows the LSTM-based framework. 
One benefit of the RL framework is that it handles non-differential cases. For example, \citet{zheng2022symbolic} proposes to learn an RNN that outputs a symbolic gradient update formula. The loss is not differentiable w.r.t. the formula, and thus RL is the only choice.

\begin{table}
    \centering 
    \scalebox{0.93}{
        \begin{tabular}{c | p{2in}|p{4in}}
        \toprule
        Framework & Paper  & Methodology \\
        \midrule
        \multirow{7}{*}{LSTM-based} & \cite{andrychowicz2016learning}  & First propose a basic LSTM-based framework.\\
                                    & \cite{lv2017learning}            & Use random scaling and convex regularizers to stabilize the training process.\\
                                    & \cite{chen2020training}          & Use curriculum learning to gradually increase the maximal number of iterations.\\
                                    & \cite{chen2020learning}          & Learn a variational stopping policy to determine the stopping time.\\
                                    & \cite{metz2019understanding}     & Use the weighted average of two unbiased gradient estimators. \\
                                    & \cite{wichrowska2017learned}     & Enhance the RNN structure in a hierarchical manner. \\
        \midrule
        \multirow{3}{*}{RL-based}   & \cite{li2016learning}            & First propose a basic RL-based framework.\\
                                    & \cite{li2017learning}           & Combine the coordinates with permutation invariance property into one
 group. \\
                                    & \cite{zheng2022symbolic} & Learn an agent that outputs a symbolic gradient update formula. \\ 
        \bottomrule 
        \end{tabular}}
    \caption{Summary of works using AI techniques to enhance gradient-based methods.} \label{tab:gradient} 
\end{table}

\subsection{Enhancement for ADMM-type methods}\label{sec:admm}
The Alternating Direction Method of Multipliers (ADMM) is an efficient
first-order optimization algorithm  that solves problems in the form
\begin{equation} \label{prob:admm}
    \min_{x \in\mR^n, s \in \mR^m}\enspace f(x) + g(s) \enspace\text{s.t.}\enspace Ax + s = b,
\end{equation}
where $A \in \mR^{m\times n}$ and $b \in \mR^m$ \citep{boyd2011distributed}. As a key component of the ADMM method, the augmented Lagrangian function for problem~\eqref{prob:admm} is given by 
\begin{equation}
    L_\rho(x, s, y) = f(x) + g(s) + y^T(Ax + s - b) + \frac{\rho}{2}\|Ax + s - b\|^2,
\end{equation} 
where $y\in\mR^m$ is the dual variable or Lagrange multiplier and $\rho > 0$ is called the penalty parameter. Then the ADMM consists of the iterations 
\begin{equation}
    \begin{split}
        x^{t+1} &= \argmin_{x \in \mR^n} L_\rho(x, s^t, y^t) \\
        s^{t+1} &= \argmin_{s \in \mR^m} L_\rho(x^{t+1}, s, y^t) \\
        y^{t+1} &= y^t + \rho(Ax^{t+1} + s^{t+1} - b).
    \end{split}
\end{equation}

Computational experiments on different applications~\citep{fortin2000augmented,fukushima1992application,kontogiorgis1998variable} have shown that, if the fixed
penalty $\rho$ is chosen too small or too large, the solution time can increase significantly. As a remedy, some heuristics have been developed to adapt $\rho$ during the optimization process, where the key idea is to balance primal and dual residuals~\citep{he2000alternating,wang2001decomposition,boyd2011distributed}. A simplified version of their adaption strategy can be summarized as 
\begin{equation} \label{eq:adaption}
    \rho^{t+1} =
\begin{cases}
     \tau \rho^t & \text{if } \|r_p^t\| \geq \mu \|r_d^t\| \\
     \frac{1}{\tau} \rho^t & \text{if } \|r_d^t\| \geq \mu \|r_p^t\| \\
     \rho^t & \text{otherwise,}
\end{cases}
\end{equation}
where 
\[
    r_p^t \coloneqq Ax^t + s^t - b \in \mR^m \enspace\text{and}\enspace r_d^t \coloneqq \rho^tA^T(s^t - s^{t-1}) \in \mR^n
\]
denote the primal and dual residuals, respectively, and $\mu, \tau > 1$ are hyperparameters.  

In this section, we want to review the literatre aiming to learn a parameterized adaption rule. \citet{zeng2022reinforcement} propose a reinforcement-learning-based penalty-tuning strategy for solving distributed optimal power flow problems \citep{mhanna2018adaptive}. 
\citet{zeng2022reinforcement} propose to learn a good penalty-tuning strategy using reinforcement learning. Following the general introduction to reinforcement learning (see Section~\ref{sec:rl}), we only need to specify the state space $\Sscr$, the action space $\Ascr$ and the reward function $R: \Sscr \times \Ascr \to \mR$. \cite{zeng2022reinforcement} set the state to include the past $k$-point history of primal and dual residuals, i.e., 
\[ s^t = \left[ (r_p^{t-k+1}, r_d^{t-k+1}), \dots, (r_p^t, r_d^t)\right] \in \mR^{k\times (m + n)} = \Sscr.\]
The problems considered by the authors have the same dimension, and thus $k, m, n$ are fixed. 
For the action space, the authors set it to be a discrete set of pre-determined possible penalty values, i.e.,
\[\Ascr = \{\rho_1, \dots, \rho_d\} \enspace\text{with}\enspace \rho_1 < \dots < \rho_d.\]
For the reward function, the authors consider two parts: termination and comparison with the baseline. For termination, the authors design the following reward function 
\begin{equation} \label{eq:termination_reward}
   R_{\text{termination}}(s^t, \rho^{t}) = 
\begin{cases}
    200 & \text{if } \|r_p^{t+1}\| < \epsilon_p \enspace\text{and}\enspace \|r_d^{t+1}\| < \epsilon_d \\
    0 & \text{otherwise,}
\end{cases} 
\end{equation}
where $\epsilon_p$ and $\epsilon_d$ are the small tolerance for checking the termination. For the comparison with the baseline, the authors set the baseline to be the regular parameter-tuning rule, e.g.~\eqref{eq:adaption}, and then they design the following reward function as the relative advantage of the RL policy over the baseline
\[
R_{\text{comparasion}}(s^t, \rho^{t}) = 
\frac{\|\hat r_p^{t+1}\| - \|r_p^{t+1}\|}{\|\hat r_p^{t+1}\|} + \frac{\|\hat r_d^{t+1}\| - \|r_d^{t+1}\|}{\|\hat r_d^{t+1}\|}, 
\]
where $\hat r_p^{t+1}$ and $\hat r_d^{t+1}$ are primal and dual residuals obtained by the baseline. Then the overall reward function is defined as the summation
\[R(s^t, \rho^{t}) = R_{\text{termination}}(s^t, \rho^{t}) + R_{\text{comparasion}}(s^t, \rho^{t}).\]

\citet{ichnowski2021accelerating} propose a reinforcement-learning-based penalty-tuning strategy for solving quadratic programming problems in the following form:
\begin{equation} \label{eq:quadratic}
     \min_{x \in \mR^n, s \in \mR^m} \frac{1}{2}x^TQx + q^T x \enspace\text{s.t.}\enspace Ax = s, \ l \leq s \leq u,
\end{equation}
where $Q\in\mR^{n\times n}$ is a symmetric positive semi-definite matrix, $q\in\mR^n$, $A\in\mR^{m\times n}$, and $l, u \in \mR^m$. 
The ADMM is one of the widely adopted optimization algorithms for solving QP Problem~\eqref{eq:quadratic}, e.g., implemented in the OSQP solver~\citep{stellato2020osqp}. Unlike the standard ADMM, the variant used in OSQP has $m$ penalty parameters, i.e., $\rho\in\mR^m$. To tackle this issue, \citet{ichnowski2021accelerating} adopted the multi-agent single-policy RL~\citep{huang2020one}. More specifically, the authors design $m$ agents for predicting entries in $\rho$ and all the agents share the same policy. For agent $i \in [m]$, the state is defined by 
\[s_i^t =
\begin{bmatrix}
    \min(s_i^t - l_i, u_i - s_i^t) \\
    (Ax^t)_i - s_i^t \\
    y_i^t \\
    \rho_i^t \\
    \|r_p^t\| \\
    \|r_d^t\|
\end{bmatrix} \in \mR^6.
\]
Thus, the state space $ \Sscr$ contains 6 elements, i.e., $ \Sscr = \mR^6$.
Note that the dimension of the state space is independent of the problem size, and thus this approach can handle different QP problems. The design of the action space is similar to~\citet{zeng2022reinforcement}. For the reward function, the authors simply set it as the termination reward as in Equation \eqref{eq:termination_reward}. 

A drawback of the approach proposed by \citet{ichnowski2021accelerating} is that the state representation is local for each agent and penalty parameter, and
thus is insufficient to capture the contextual information. To resolve this issue, \citet{jung2022learning} extended this framework by using the graph representation of QP as the state representation and then the corresponding $Q$ function is modelled by message-passing on the graph. See Section~\ref{sec:gnn} for a more detailed introduction to the graph representation of QP. 

\begin{table}
    \centering 
        \begin{tabular}{c | l | l}
        \toprule
        Paper & Problem type  &  Methodology \\
        \midrule
        \cite{zeng2022reinforcement} & distributed OPF & RL\\
        \cite{ichnowski2021accelerating} & QP & Multi-agent single-policy RL\\
        \cite{jung2022learning} & QP & GNN + RL\\
        \bottomrule 
        \end{tabular}
    \caption{Summary of works using AI techniques to enhance ADMM-based methods.} \label{tab:admm} 
\end{table}

\subsection{Enhancement for column-generation methods} \label{sec:column_generation}

\begin{figure}
    \centering
    \includegraphics[width=\linewidth]{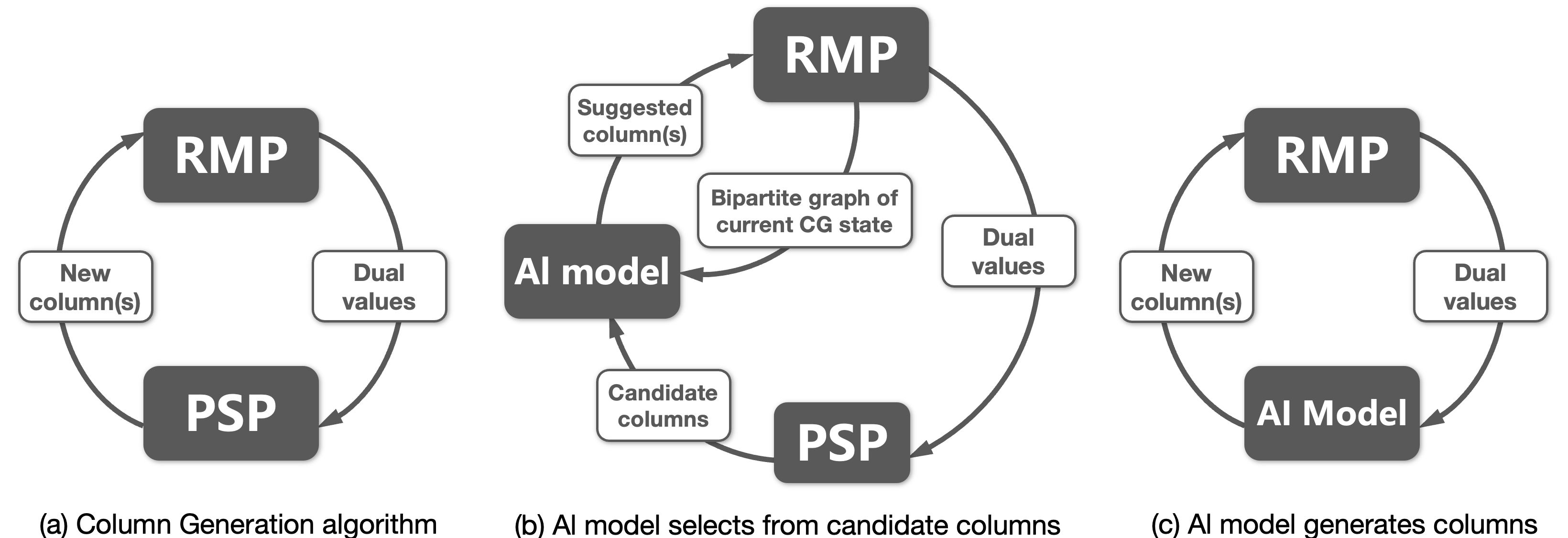}
    \caption{Illustration for column generation (CG) algorithm and how does AI model assist CG. The credits for sources (a) and (b) are attributed to \citet{chi2022deep}. }
    \label{fig:CG-AI}
\end{figure}

Column generation (CG) is an algorithm for solving linear programs (LPs) with a prohibitively large number of variables (i.e., columns)~\citep{desaulniers2006column}. 
CG starts by solving a restricted master problem (RMP) with a subset of columns and gradually generates new columns that can improve the solution of the current RMP. The method stops when no such columns exist. More specifically, we consider the following master LP problem
\begin{alignat*}{6}
 \text{[MP] } \quad   \min_{x\in\mR^{|\Omega|}} \enspace & \enspace \sum_{p\in\Omega}c_px_p \\
    \text{s.t.} \enspace & \enspace \sum_{p\in\Omega}x_p\mathbf{a_p}=\mathbf{b}\\
                         & \enspace x_p \geq 0, \enspace\forall p \in \Omega,
\end{alignat*}
where $\Omega$ represents the set of variable indices. We consider the case when $|\Omega|$ is large and the variables can not be enumerated explicitly. The CG algorithm starts with a subset $\Fscr \subset \Omega$ of variables and deals with the following restricted master problem [RMP] 
\begin{alignat*}{6}
\text{[RMP] } \quad    \min_{x\in\mR^{|\Fscr|}} \enspace & \enspace \sum_{p\in\Fscr}c_px_p \\
    \text{s.t.} \enspace & \enspace \sum_{p\in\Fscr}x_p\mathbf{a_p}=\mathbf{b} \\
                         & \enspace x_p \geq 0, \enspace\forall p \in \Fscr.
\end{alignat*}
Let $y_\Fscr \in \mR^m$ denote the optimal variable for the [RMP]. 
Then new columns with negative reduced cost will be generated by solving another subproblem, called the pricing subproblem (PSP) or column generation subproblem (CGSP). 
After solving [PSP], new column(s) $\Gscr$ is generated. $\Gscr$ must be negative-reduced-cost columns, i.e., $
    \Gscr \subseteq \{ p\in\Omega \mid  c_p - y_\Fscr^T\mathbf{a_p} < 0\}
$. 
In practice, enumerating all possible columns $p\in\Omega$ is computationally infeasible. 
Thus, ${a_p}$ is usually implicitly specified by the constraint of pricing problem [PSP]. 
The objective of the [PSP] is to incorporate dual values and identify the most promising column(s) that can augment the [RMP] solution. 
Then, the CG algorithm updates the index set by $\Fscr \leftarrow \Fscr \cup \Gscr$, and repeats the above process until no such column exists. 
The bottleneck of CG is [PSP], since [PSP] usually involves integer variables. Meanwhile, the total number of CG iterations largely depends on the strategy of choosing the negative-reduced-cost columns. An illustration of the CG algorithm is shown in Figure~\ref{fig:CG-AI}(a). 

\citet{Desaulniers2020MachineLearningBasedCS} pioneered the research direction on using AI techniques to {select} columns in the CG algorithm. 
As shown in Figure~\ref{fig:CG-AI}(b), the [PSP] first heuristically or exactly provides candidate columns $\bar{\Gscr} \subseteq \{ p\in\Omega \mid  c_p - y_\Fscr^T \mathbf{a_p} < 0\}$. For optimization problems like the vehicle and crew scheduling problem and the vehicle routing problem with time windows, it is possible to generate {multiple} negative-reduced-cost columns at no extra computational effort. 
Given the candidate columns $\overline{\Gscr},$
we want to select a proper subset $\Gscr \subseteq \overline{\Gscr}$. The challenge in this task is that we don't want the size of $\Gscr$ to be very large as it may make the next restricted master problem hard to solve. In the meantime, we want $\Gscr$ to contain more ``promising" columns. \citet{Desaulniers2020MachineLearningBasedCS} suggest that the importance of columns can be obtained by solving the following mixed integer linear program (MILP):
\begin{subequations} \label{eq:milp-pp}
\begin{align} 
    \min_{x\in\mR^{|\Fscr\cup\overline{\Gscr}|}, y \in \{0,1\}^{\overline{\Gscr}}} \label{eq:milp-pp-obj}
    \enspace &
    \enspace \sum_{p\in\Fscr \cup \overline{\Gscr}}c_px_p + \lambda \sum_{p\in\overline{\Gscr}}y_p \\
    \text{s.t.} \enspace & \enspace \sum_{p\in\Fscr \cup \overline{\Gscr}}x_p\mathbf{a_p} = \mathbf{b}, \\
                         & \enspace x_p \geq 0, \enspace\forall p \in \Fscr \cup \overline{\Gscr}, \\
                         & \enspace x_p \leq y_p, \enspace\forall p \in \overline{\Gscr}, \label{eq:milp-pp-swither}
\end{align}
\end{subequations} 
where the binary variables $y_p$ decide whether the column $p$ should be selected for the next RMP. 
More specifically, after solving problem~\eqref{eq:milp-pp}, if $y_p = 1$, then column $p$ is a candidate for selection, otherwise, column $p$ is excluded.
The hyperparameter $\lambda$ in \eqref{eq:milp-pp-obj} balances the trade-off between minimizing the [RMP] objective and controlling the size of $\Gscr$. 
In other words, by solving Problem \eqref{eq:milp-pp}, we can derive a better subset $\Gscr$ with two advantages. Firstly, the subsequent [RMP] objective is greatly reduced, leading to an improved solution for the master problem. Secondly, the size of the resulting subset remains manageable.
Once the MILP is solved, the better subset $\Gscr$ is determined by:
\[
\Gscr = \{p \in \overline{\Gscr} \mid y_p = 1\}.
\]

However, solving a MILP at each iteration can be time-consuming in practice. Thus the authors suggest imitating this MILP using AI techniques. More specifically, they represent the MILP as a bipartite graph as the one introduced in Section~\ref{sec:gnn} (seen in Figure \ref{fig:lp-graph}). 
Then they train a standard GNN to imitate the optimal $\{y_p\}_{p\in \bar{\Gscr}}$ from the MILP~\eqref{eq:milp-pp}. 
This is a binary classification task. In the inference stage, the GNN will predict the probability that a generated column $p$ should be selected, i.e., $Pr(y_p=1)$. If this probability is greater than $0.5$, then column $p$ is considered promising and added to the [RMP] in the next iteration.

Alternatively, \citet{chi2022deep} propose a reinforcement-learning-based approach to select columns in the CG algorithm. 
More specifically, at each iteration of the CG algorithm, the {state} $s$ is the bipartite graph representing the current [RMP] plus current iteration information, including the candidate columns $\overline{G}$ with their reduced costs and solution values, and how long a column stays or leaves the basis. 
In each iteration, the authors want to select {one} column from the candidate columns, so the action space $\Ascr$ is  $\overline{G}$. 
To encourage the RL agent to select columns such that the CG algorithm converges {better} and {faster}, the {reward} function at each step consists of two components: 
{1)} The first component $\frac{obj_{t-1} - obj_t}{obj_0}$ gives a higher score if taking action $a_t$ leads to a {better} outcome, i.e., it causes a faster decrease in the objective value. Here, $obj_t$ is the objective value of the RMP at time step $t$, and $obj_0$ is the objective value of the [RMP] in the first CG iteration. The latter is used to normalize $(obj_{t-1} - obj_t)$, ensuring that the first component is numerically comparable across different LP instances. 
{2)} The second component is simply $-1$, which encourages the agent to prefer {shorter} iterations. 
This component represents the total number of iterations in the cumulative rewards. 
To balance these two components, a non-negative hyperparameter $\alpha$ is introduced. The reward at time step $t$ is then defined as 
\[R_t = \alpha\frac{obj_{t-1} - obj_t}{obj_0} - 1.\]

\citet{chi2022deep} also use the bipartite graph to represent the [RMP]. Given the bipartite graphs with node features, GNNs are used as the Q-function approximation. 
Given a particular state, the Q-function estimates the expected future reward (Q-values) for all possible actions, i.e., all candidate columns. 
GNNs are capable of capturing the complex relationships between nodes in the graph, which makes them suitable for this task. 
Note that \citet{chi2022deep} select only {one} column at each iteration, so the action with the maximum Q value is selected to add to the [RMP]. 
Since the Q-function will consider future rewards, the agent can make a better column selection at each step. 
Compared with the heuristic (greedy) strategy for selecting columns, the RL-based algorithm converges faster in terms of the number of iterations and total time.   
\citet{chi2022deep} also mention that adopting a curriculum learning paradigm improves the learning of the RL agent and results in better column selection and faster convergence. This is crucial to convergence when the training set contains instances of varying difficulties.

Besides selecting from candidate columns, the AI model can also directly generate columns. 
\citet{shen2022enhancing} consider the graph coloring problem. A column corresponds to a maximal independent set (MIS). The pricing problem for generating a MIS is an NP-hard Maximum Weight Independent Set Problem. 
\citet{shen2022enhancing} leverage an AI-based primal heuristic to generate a near-optimal MIS from the pricing problem.
Given the feature vector of each vertex, this AI-based heuristic estimates the probability of each vertex being selected into the MIS. It is then used to guide a sampling method to efficiently generate multiple high-quality MISs (i.e., columns).
The sampling method ensures the column is indeed an MIS by sequentially adding vertices, marking adjacent vertex as invalid, and adjusting the probability of selecting remaining vertices. 
The diversity of multiple MISs is achieved by randomly selecting the starting vertex. 
Due to generating high-quality columns efficiently, MLPH significantly accelerates the progress of CG and reduces CG iterations, especially for larger and denser graphs.  

To conclude, the different AI models assist CG by either selecting more promising columns or directly generating columns. 
The former category naturally ensures the candidate columns are valid and can reduce the CG iterations further if RL or curriculum learning is adopted. However, the time for solving [PSP] is irreducible. 
In contrast, the latter category replaces solving [PSP] with an AI-model inference and thus has greater potential. 
If the latter category is further augmented with an RL algorithm, it may achieve greater end-to-end speedup. 
Its drawback is that the AI model cannot directly generate a column feasible to [PSP]. An ad-hoc sampling step has to be designed to mitigate the gap between AI prediction and feasible column generation, for each type of problem and the corresponding [PSP]. 
It would be interesting to explore the possibility of either designing a unified sampling step or developing an end-to-end AI model capable of generating the columns. 

\subsection{Enhancement for simplex method} \label{sec:simplex}
The simplex algorithm is a classical method for solving linear programming problems of the form 
\begin{align*}
    \min_{x\in\mR^n} \enspace & \enspace c^Tx \\
    \text{s.t.} \enspace & \enspace Ax = b, \\
                         & \enspace x \geq 0.
\end{align*}

In the following, we will first briefly introduce the primal simplex algorithm.
It consists of two crucial components: basis initialization and pivoting. 
Typically, the introductory textbooks assume that the primal simplex algorithm begins with a primal feasible basis. 
This assumption is reasonable as techniques like the big-M method are able to address the infeasibility. 
Nonetheless, it remains a question whether the initial basis is ``good," i.e., whether it is close to and converges rapidly to the optimal solution. 
After basis initialization, the primal simplex algorithm repeats a pivoting process, selecting the entering basis variable and the leaving basis variable. 
Geometrically, each basis corresponds to a feasible solution and a vertex in the feasible region (a polyhedron). 
Accordingly, the simplex algorithm starts from the initial solution vertex and each pivot corresponds to a step transition towards adjacent vertices until it reaches the optimal solution vertex. 
To be specific, 
The simplex method starts from an initial feasible basis $\Bscr=[\Bscr_x; \Bscr_s]$ containing the indices of basic variable and constraint variables. 
At each iteration, we first divide the variable $x \in \mR^n$ (assume the slacks are also included) into $x = [x_B; x_N]$, where $x_B$ and $x_N$ represent the basic and non-basic variables, respectively. The matrix $A = [B; N]$ and the vector $c = [c_B; c_N]$ are also divided correspondingly, where $B$ is the matrix of coefficients for the basic variables and $N$ is the matrix of coefficients for the non-basic variables. It is required that each column of the matrix $B$ be linearly independent. Substituting the constraint $Ax = b$ into the objective function, we obtain $c^T x = c_B^T B^{-1} b + (c_N^T - c_B^T B^{-1}N)x_N$, where $c_B^T$ and $c_N^T$ represent the transpose of the vectors $c_B$ and $c_N$, respectively. The second term is called the reduced cost vector $\bar{c} = c_N^T - c_B^T B^{-1}N$, which represents the cost of increasing the value of the non-basic variables. The non-basic variables corresponding to the negative components of $\bar{c}$ can cause a decrease in the objective function. Therefore, the selection of the entering basis variable is the process of selecting a non-basic variable corresponding to a negative component of $\bar{c}$. Different pivot rules provide different methods for selecting the entering basis variable. 
The essence of the pivot rule of the simplex algorithm is to convert the status of a certain column between basic and non-basic. 
Once the basic variables are determined, we can use $x_B = B^{-1} b - B^{-1}Nx_N \geq 0$ to derive the leaving basis variable. The pivoting process is repeated until the basis corresponding to the optimal solution is obtained, i.e., the end of the simplex algorithm. 

For basis initialization, \citet{fan2023smart} propose a GNN-based initial basis selection strategy. 
They first represent an LP as a bipartite graph and convert the basis selection task into a classification task, i.e., determine the basis status for each constraint and variable. 
To ensure the predicted basis status always satisfies the bound requirements, e.g., a free variable (with unbounded lower and upper bounds) can only be basic and a variable with an unbounded upper bound will not achieve its upper bound. 
Then in the inference stage, \citet{fan2023smart} propose the basis generation and adjustment steps ensuring the basis is valid, i.e., the corresponding constraint matrix is non-singular. 
Extensive experiments are conducted, including large-scale industrial cases, 
to demonstrate the performance of the GNN-based initial basis selection strategy over traditional initialization heuristics that fail to produce a close-to-optimal basis  \citep{ploskas2021triangulation}. In contrast, the GNN-based method utilizes past solved linear programs and smartly builds a close-to-optimal basis.
 
The selection of the entering basis variable, commonly known as the pivoting strategy, plays a crucial role in determining the efficiency of the simplex method. Various classical pivoting strategies \citep{dantzig2003linear,forrest1992steepest} have been proposed in the literature, and their effectiveness has been evaluated empirically. The Dantzig pivoting rule \citep{dantzig2003linear} is a popular strategy that selects the non-basic variable with the most negative reduced cost as the entering variable. In contrast, the steepest edge pivoting rule \citep{forrest1992steepest} selects the non-basic variable with the largest rate of decrease of the objective value per unit distance travelled along the improving edge. The choice of pivoting strategy is highly dependent on the specific problem at hand, and different strategies may yield significantly different results. \citet{suriyanarayana2022reinforcement} use reinforcement learning techniques to switch between Dantzig’s rule and the steepest edge rule. Namely, at each iteration of the simplex algorithm, the trained agent will select one of the two above-mentioned pivoting rules. Alternatively, \citet{li2022rethinking} aims to learn a new pivoting strategy through the application of the Monte Carlo Tree Search (MCTS) method. Their contribution focuses on four core aspects, including transforming the simplex method into a pseudo-tree structure, constructing appropriate reinforcement learning models, finding the optimal pivot sequence under the guarantee of theory, and providing a complete method for discovering multiple optimal pivot paths. The study proposes a novel imitative tree structure, SimplexPseudoTree, for the exploration of optimal pivot paths, and constructs four reinforcement-learning models to determine the optimal pivot paths based on the MCTS method. The research provides theoretical as well as computational experiments to demonstrate the optimality of the proposed MCTS rule. 
The MCTS rule can avoid unnecessary searches and determine the shortest pivot paths of the simplex method, leading to more efficient problem-solving in the context of linear programming, namely solving the problem with fewer iterations.

\section{Algorithm Selection and Design for Discrete Optimization} \label{sec:branch_bound}
Consider an optimization problem of minimizing $f$ over some finite set $\Xscr$, i.e., 
\[
\min_{x\in\Xscr} \enspace f(x).
\]
The branch-and-bound (B\&B) method, initially proposed by \citet{land2010automatic}, recursively divides the finite feasible region $\Xscr$ into its subsets $\Xscr_1, \Xscr_2, \dots, \Xscr_p$ until no more division is possible such that 
\[\Xscr = \bigcup_{i=1}^p \Xscr_i.\]
All these subsets form a B\&B tree. The key assumption in the B\&B method is that for every subset or node $\Xscr_i$, there is an algorithm that can calculate 
\begin{itemize}
    \item A lower bound, i.e.,  
    $\underline{f_{\Xscr_i}} \leq \min_{x\in\Xscr_i} f(x).$
    \item An upper bound provided by each feasible point $\overline{x} \in \Xscr_i$, i.e., 
    $\overline{f_{\Xscr_i}} = f(\overline{x}) \geq  \min_{x\in\Xscr_i} f(x).$
\end{itemize}
The basic idea behind the B\&B algorithm is that, if for two subsets $\Xscr_1, \Xscr_2 \subseteq \Xscr$,
\[\overline{f_{\Xscr_2}} \leq \underline{f_{\Xscr_1}},\]
then the solutions in $\Xscr_1$ can be disregarded. More specifically, in each iteration, the B\&B algorithm will pick a child node $\Xscr_i$, which is called node selection, and check whether the lower bound exceeds the best available upper bound. If it exceeds, then the B\&B algorithm can safely remove the current node from the B\&B tree, and continue with the next child node. If it doesn't exceed, then the B\&B will divide the current node by branching on a variable. This step is also called variable selection. 

\subsection{Mixed-integer linear programming}

\begin{figure}
    \centering
    \includegraphics[width=\linewidth]{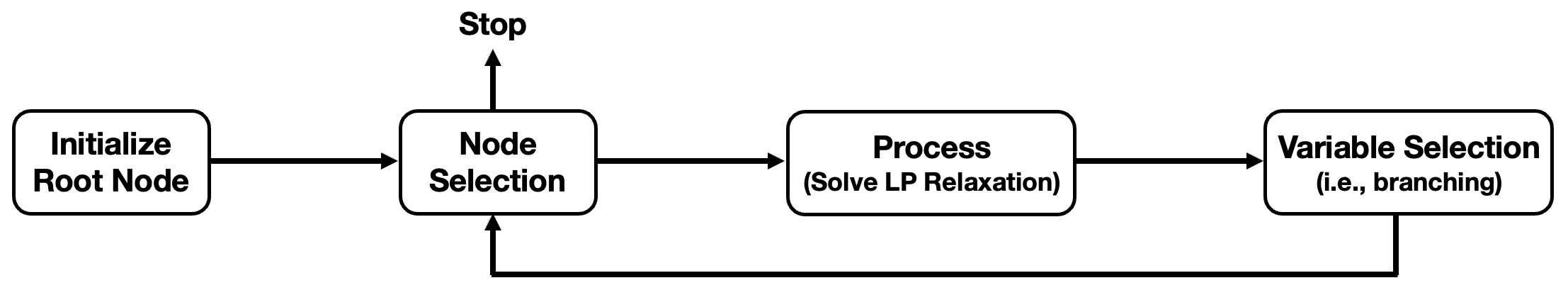}
    \caption{
    The steps of the branch-and-bound algorithm for mixed-integer linear programming problems.
    }
    \label{fig:branch_and_bound}
\end{figure}

Mixed-integer linear programming problems are the core of discrete optimization because they can model a wide variety of problems in different applications. A MIP is an optimization problem of the form
\begin{subequations} \label{eq:mip}
\begin{align} 
    p^* = \min_{x\in\mZ^d \times \mR^{n-d}} \enspace & \enspace c^Tx \\
    \text{s.t.} \enspace & Ax = b, \label{eq:mip-b} \\
                         & x \geq 0, \label{eq:mip-c}
\end{align}
\end{subequations} 
where $A \in \mR^{m \times n}$. When applying the B\&B algorithm to solve a MIP problem, one first relaxes the integrality constraint and obtain a linear program (LP)
\begin{subequations} \label{eq:mip-relax}
\begin{align}
    \overline{p} = \min_{x\in\mR^{n}} \enspace & \enspace c^Tx \label{eq:mip-relax-a} \\
    \text{s.t.} \enspace & Ax = b, \label{eq:mip-relax-b} \\
                         & x \geq 0, \label{eq:mip-relax-c}
\end{align}
\end{subequations}
whose solution provides a lower-bound, i.e., $\overline{p} \leq p^*$. Let $\overline{x}$ denote the minimizer for the LP relaxation. If $\overline{x}$ satisfies the integrality constraint for the original MILP problem, then $\overline{x}$ is a global solution, and we are done. Otherwise, we need to decompose the LP relaxation into two sub-problems by selecting a variable violating the integrality constraint, i.e., 
\[j \in [p] \enspace\text{such that}\enspace \overline{x}_j \notin \mZ.\]
The two sub-problems have the form
\begin{alignat*}{6}
    p_{j}^- = \min_{x\in\mR^{n}} \enspace & \enspace c^Tx \\
    \text{s.t.} \enspace & Ax = b, \\
                         & x_i \geq 0 \enspace \forall i \neq j, \\
                         & 0   \leq \ x_j \leq \floor*{\overline{x}_j}
\end{alignat*}
and
\begin{alignat*}{6}
    p_{j}^+ = \min_{x\in\mR^{n}} \enspace & \enspace c^Tx \\
    \text{s.t.} \enspace & Ax = b, \\
                         & x_i \geq 0 \enspace \forall i \neq j, \\
                         & x_{\overline{i}} \geq \ceil*{\overline{x}_{j}}.
\end{alignat*}
The variable selection then refers to the process of selecting variable $j$. Then the B\&B algorithm is going to pick one of the child LP problems and continue the above process. The node selection then refers to picking the child LP problem.  

Node and variable selection largely affect the performance of the B\&B algorithm~\citep{huang2021branch}. Many machine-learning-based approaches have been developed to assist node selection~\citep{He2014,song2018learning,sabharwal2012guiding} and variable selection~\citep{Khalil2016,Alvarez2017,di2016dash,balcan2018learning,Gasse2019,Gupta2020,gupta2022lookback,zarpellon2021parameterizing,qu2022improved,etheve2020reinforcement,sun2020improving} in the B\&B algorithm for solving MIP problems. In the next two sections, we provide an overview of these approaches. 

\subsubsection{Variable selection} 
Variable selection is a crucial task in the branch-and-bound algorithm. It determines the way in which a current node is partitioned into two child nodes in a recursive manner, by choosing which fractional variables, also known as candidate variables, to branch on. The goal of an effective branching strategy is to minimize the number of explored nodes before the termination of the B\&B algorithm. The variable selection plays a vital role in the B\&B method and has led to the development of various heuristics over time~\citep{achterberg2005branching,linderoth1999computational}.

The simplest heuristic for variable selection is the {most-infeasible branching} rule~\citep{achterberg2005branching}. It suggests branching on the variable with the greatest fractional part. For binary variables, this corresponds to selecting the variable whose value in the LP relaxation is furthest from being an integer, e.g., $0.5$. 
The intuition behind this is that by branching on the most fractional variable, the most ``ambiguous'' part of the current solution {may} be prioritized, hoping to quickly converge to an integer solution.
However, this method has been shown to perform poorly in practice. Another popular heuristic is {pseudocost branching}~\citep{benichou1971experiments}, which uses a history of increase in the dual bounds observed during previous branching to estimate the dual bound improvements for each candidate variable at the current node. Although its performance improves as the B\&B algorithm progresses, it usually performs poorly in the early stages. {Strong branching}~\citep{applegate1995finding} is another well-known heuristic. It evaluates the dual bound increase for each fractional variable by computing the linear programming relaxations resulting from branching on that variable. The variable that results in the largest increase is then selected as the branching variable for the current node. Although strong branching can produce a B\&B tree with a small number of nodes, its high computational cost often makes it intractable in practice.

In recent years, several machine-learning-based variable selection strategies have been developed. These strategies can be divided into three categories: 1) models that switch between different branching rules~\citep{di2016dash,balcan2018learning}, 2) models that mimic a strong but expensive branching rule~\citep{Khalil2016,Alvarez2017,Gasse2019,Gupta2020,gupta2022lookback,zarpellon2021parameterizing}, and 3) the utilization of reinforcement learning to learn a new branching strategy~\citep{qu2022improved,sun2020improving,etheve2020reinforcement}. These AI-based approaches aim to overcome the limitations of the classical variable selection rules, such as the most-infeasible branching and pseudocost branching, by providing a more adaptive and flexible approach to variable selection in the B\&B algorithm.

\paragraph{Switching Between Different Branching Rules} 
\citet{di2016dash} conducted a study on the dynamic and sequential nature of branch-and-bound algorithms used in mixed-integer linear programming problems. They show that no single branching rule could perform optimally across different subproblems of the same MILP. This observation motivated the development of the Dynamic Approach for Switching Heuristics (DASH) algorithm. DASH employs a two-step approach, where the first step involves clustering problems based on defined features using the K-means algorithm, and the second step involves learning the correct assignment of branching rules to each cluster during an offline training phase. The algorithm adapts to changes in the instance as the search depth increases, switching to a new branching rule that best fits the current cluster. 

Rather than selecting only one branching rule, an alternative approach proposed by \citet{balcan2018learning} is combining multiple branching rules that are score-based.
It means the rule relies on a quantitative ``score" assigned to each variable to determine its priority in the branching process. 
Take the most-infeasible branching rule mentioned above as an example. 
Denote the solution of the LP relaxation associated with the current node as $x$. 
The score of $i$-th candidate variable is $\text{score}(x_i) = \min(\lceil x_i \rceil - x_i, x_i - \lfloor x_i \rfloor)$. 
Given the scores from different score-based branch rules, 
\citet{balcan2018learning} proposed a learning-based approach where scores are combined and weighted to enhance accelerate the B\&B algorithm. 


\paragraph{Imitation Learning} 
As mentioned earlier, strong branching is a well-known variable selection strategy in B\&B algorithms, which is known for its remarkable practical performance. More specifically, let 
\[\Cscr \subseteq \{i \in [p] \mid \overline{x}_i \notin \mZ\}\]
denote the set of branching candidates. The key idea of strong branching is to calculate a score $s_i$ for every possible candidate $i \in \Cscr$ and then select the one with the highest score, i.e., 
\[\overline{i} = \argmax_{i \in \Cscr} \,\, \text{score}_i.\]
Typically, the score is defined by the improvement in the lower bound, i.e., 
\[\text{score}_i = \max\{p_i^- - \overline{p}, \epsilon\}\cdot\max\{p_i^+ - \overline{p}, \epsilon\},\]
where $\epsilon$ is some small constant. However, the implementation of strong branching is computationally demanding, as it requires the resolution of two LP problems for each candidate variable. To address this issue, several researchers have explored the use of AI techniques to mimic the strong branching rule in B\&B algorithms. \citet{Khalil2016} proposed the first work in this direction, where they developed an SVM model that learns a branching rule customized to a single instance during the B\&B process. This process involves collecting a set of B\&B nodes and performing strong branching on these nodes to obtain the ranking of the candidate variables. The variables are then categorized into ``good" and ``bad" based on the scores, and the SVM model is trained to identify ``good" variables. \citet{Alvarez2017} introduced a two-phased approach, which results in a ``learned" branching strategy that can be used as an approximation of strong branching within the B\&B algorithm. The first phase involves solving a set of training problems with strong branching as a branching heuristic and recording each branching decision in a training set. In the second phase, the learned heuristic is introduced into B\&B and evaluated for efficiency on a set of test problems. 
\citet{Gasse2019} utilized a GNN model to tackle the variable selection problem in B\&B. This model takes in the bipartite graph representation of an MILP. 
This representation includes nodes representing both constraints and variables, along with their respective features and connectivity. 
To be specific, constraint node features can be the dual solutions from LP relaxation, cosine similarities between each constraint's coefficients and objective coefficients, 
The variable node features can be the objective coefficients and lower and upper bounds of variables. 
Due to the graph representation and the GNN model, MILPs with varied sizes can be handled. 
The GNN is trained to approximate strong branching using imitation learning and has been shown to improve upon previous branching approaches for several MILP problem benchmarks and is competitive with state-of-the-art B\&B solvers. 
However, this method requires a high-end GPU card to speed up the GNN inference time, which is not always feasible for MILP practitioners. To overcome this limitation, \citet{Gupta2020} studied the time-accuracy trade-off in learning to branch, and proposed a hybrid architecture that uses a GNN model at the root node and a fast but weak predictor, such as a Multi-Layer Perceptron (MLP), at the remaining nodes. This approach enhances the weak model with high-level structural information extracted at the root node by the GNN model. \citet{gupta2022lookback} later found that the strong branching heuristic often leads to a child node's best choice being the parent's second-best choice, known as the ``lookback" phenomenon. To imitate this behaviour more closely, they proposed two methods that incorporate the lookback phenomenon into GNN training by adding a lookback regularized term. 
Finally, \citet{zarpellon2021parameterizing} also inherit the imitation learning framework, but innovatively propose to utilize the information of B\&B search trees. 
They believe that many MILPs share similarities in terms of the search trees. 
However, there is no natural parameterization of the search tree. 
To address this gap, a set of 61 hand-crafted input features is proposed to describe candidate variables in terms of their roles in the B\&B process. 
These tree features capture various aspects, such as the current node's depth and bound quality, the tree's growth rate and composition, the evolution of global bounds, aggregated variables’ scores, statistics on feasible solutions, and depths of open nodes. 
Experimental evidence suggests that explicitly incorporating these features enhances the accuracy of the learned agent and the agent effectively helps reduce the size of the search tree. 
The proposed method outperforms the current state-of-the-art approach and allows for generalization to unseen MILP instances. This generalization empirically verifies that unseen instances share similar search trees to training instances. 

In conclusion, all of the above works aim to mimic the strong branching rule for solving MIPs from a specific domain and have achieved promising results. However, more research is needed to fully understand the potential and limitations of these methods to advance the field of variable selection in B\&B algorithms.

\paragraph{Reinforcement Learning} 
As stated in Section \ref{sec:rl}, one limitation of imitation learning is that the performance of the learned strategy is limited by the expertise of the expert. In the following, We will first describe how the expert (i.e., strong branching) is limited and then discuss the reinforcement learning-based approach proposed by \citet{sun2020improving}. 

Strong branching is commonly acknowledged as an effective algorithm, primarily because it often results in the smallest search trees in the B\&B algorithm. This compactness in search trees can be attributed to two main factors: decision quality and the impact of the decision on other nodes. 
{1) Decision Quality}: strong branching involves making decisions on variable selection. These decisions are considered of high quality, as they contribute to faster convergence towards the optimal solution. 
{2) Impact on other nodes}: at each node in the search, strong branching will solve several relaxed LPs before making decisions on variable selection. In the process of solving relaxed LPs, relevant information is produced. This relevant information, by default, is not discarded and helps accelerate B\&B algorithm. We call the utilization of relevant information as secondary effects. 
For example, a) strong branching can evaluate the pruning conditions prior to actual branching and relevant information can help prune some subproblems. b) When solving the LP of the current branch, relevant information is obtained that can enhance other LP relaxations by eliminating unnecessary constraints. 

Among these two factors, we may assume that decision quality plays a more significant role in accelerating B\&B algorithm. However, the empirical findings of \citet{sun2020improving} suggest otherwise. 
In their study, the authors disabled the secondary effects in the full strong branching algorithm and observed that the reduction in tree size was notably less substantial compared to when the secondary effects were enabled. 
This implies that the acceleration brought by strong branching primarily stems from the secondary effects. 
Meanwhile, imitation learning can only imitate the decisions made by strong branching, not the secondary effects. 
As a result, imitating strong branching may not be a wise choice for learning a variable selection policy.

In response to these findings, \citet{sun2020improving} proposed a reinforcement learning-based approach to model the variable selection process as a Markov Decision Process (MDP). 
The authors design a primal-dual policy network, which is similar to the GNN model operating on bipartite representation. 
By setting the cumulative reward as the negative value of the total number of solving nodes within the B\&B algorithm, the learned policy is non-myopic, aiming to solve problems in fewer steps. 
Moreover, to encourage exploration during the learning process, they introduced a novelty score of the current policy. 
This score is determined by examining the policy's surrounding neighbourhood. A policy is deemed novel if it significantly deviates from its neighbours. 
The novelty score is integrated into the cumulative reward, guiding policy evolution. Such a combination enables the RL agent to navigate novel states in the B\&B process, bypass local optima, and adopt a variety of strategies.

Another limitation of imitation learning lies in the mismatch between demonstration and real data. 
Many studies employ imitation learning to emulate the strong branching method and rely solely on data gathered from expert policies. 
However, when the learned policy is applied to unseen instances, it might not always make decisions as accurate as strong branching. Consequently, the resulting states can diverge from the training data. The experts cannot provide demonstrations for every potential state the model might encounter as problem characteristics and structures can vary significantly. Thus, there might exist a mismatch between training data and unseen problem instances.
To address this, \citet{qu2022improved} proposed a novel reinforcement learning-based branching algorithm that trains on training data at the early stage to accelerate the learning process. The model then updates with a mixture of training and self-generated data to balance the exploration and exploitation. 
As a byproduct, this approach was also found to overcome the issue of the large variance in gradient estimation, which is a common challenge in MDP-based approaches. 

\begin{table}
    \centering 
        \begin{tabular}{p{0.13\linewidth} | p{0.25\linewidth} | p{0.25\linewidth} | p{0.28\linewidth} }
        \toprule
        Paper & Goal  & Methodology & Features\\
        \midrule
        \cite{di2016dash} & Learn to switch between different variable selection strategies & Cluster MILP problems into clusters using the K-means algorithm and learn an assignment of branching methods to each cluster & 40 features describing different statistics of the problem\\
        \midrule
        \cite{balcan2018learning} & Learn to combine the scores returned by different existing variable selection strategies & Empirical risk minimization to find the optimal weight for convex combinations & Similar to algorithm configuration, it finds an optimal weight for a certain distribution of MILPs. Since the weight only applies to this distribution, it requires no instance-specific features. \\
        \midrule
        \cite{Khalil2016} & Learn to mimic the strong branching strategy &  Imitation learning via SVM & 18 static features and 54 dynamic features\\
        \midrule
        \cite{Alvarez2017} & Learn to mimic the strong branching strategy & Imitation learning via random forests & static problem features, dynamic problem features and dynamic optimization features \\
        \midrule 
        \cite{Gasse2019} &  Learn to mimic the strong branching strategy & Imitation learning using GNN & Bipartite graph representation with 5 features for the constraint, 13 features for the variable and 1 feature for the constraint matrix\\
        \midrule
        \cite{Gupta2020} & Learn to mimic the strong branching strategy & Imitation learning using GNN at the root node and SVM at remaining nodes & same features as \cite{Gasse2019} at the root node and same features as \cite{Khalil2016} at remaining nodes \\
        \midrule 
        \cite{gupta2022lookback} & Learn to mimic the strong branching strategy & Imitation learning via GNN and Lookback regularization & same features as \cite{Gasse2019} \\
        \midrule
        \cite{zarpellon2021parameterizing} & Learn to mimic the SCIP's default branching strategy & Imitation learning using deep neural network & 25 features from the candidate variables and 61 features describing the state of the B\&B search tree\\
        \midrule
        \cite{sun2020improving} & Learn a novel branching strategy & Reinforcement learning with GNN &  same features as \cite{Gasse2019} \\
        \midrule
        \cite{qu2022improved} & Learn a novel branching strategy & Reinforcement learning with Double DQN &  same features as \cite{Gasse2019} \\
        \bottomrule 
        \end{tabular}
    \caption{Summary of works using AI techniques to enhance variable selection in the B\&B algorithm.} \label{tab:b&b_vs} 
\end{table}

\subsubsection{Node selection}
When applying the B\&B algorithm to solve MILP problems, the solving process involves breaking down a problem into smaller sub-problems, referred to as nodes, and selecting which node to process next. This selection is based on two main goals: finding good feasible MILP solutions to improve the upper bound and getting good LP relaxations to improve the lower bound. In the literature, various search methods have been proposed for node selection in the B\&B algorithm. One of the earliest methods is depth-first search, proposed by \citet{dakin1965tree}, where the node with the maximum depth in the B\&B search tree is selected. This method is efficient in terms of memory consumption. Another popular node selection method is the best-first search, which is proposed by \citet{hart1968formal}. This method selects the node with the currently smallest dual objective value. 

Recently, various AI-based node selection strategies have been proposed to enhance the performance of classical node selection strategies in the branch-and-bound algorithm. 
\citet{He2014} introduced an imitation learning method that learns a node selection strategy by observing a small set of solved problems. 
The method assumes that the problems at the test time exhibit similar characteristics, such as problem type, size, and parameter distribution, as those observed during the training time. The node selection policy is designed to repeatedly pick a node from the queue of unexplored nodes in a manner that mimics a simple oracle, which knows the optimal solution in advance and only expands nodes containing the optimal solution.

However, in practice, many MILPs are substantially large. This poses a significant challenge for the implementation of the imitation learning methods discussed earlier \cite{He2014}, as constructing a training set based on optimal solutions would be extremely time-consuming.
To overcome this challenge, \citet{song2018learning} presented a variant approach of imitation learning for node selection. 
In this approach, the expert will use a cut-off technique, i.e., a solver runs until a certain node limit is reached and outputs the best solution found. 
Then, the shortest path to the resulting solution becomes the expert demonstration of node selection trajectory. 
In addition to the cut-off technique, a ``gradual scaling up'' technique is applied to enhance the agent's scalability and generalization in tackling larger problem instances.  
Specifically, after the agent is trained on problems of a certain size, larger problems are generated for the agent to interact with. During this interaction, the agent provides node selection suggestions within the B\&B algorithm, and the best solution found so far serves as expert feedback to further enhance the agent's performance. 
To summarize, this particular variant of imitation learning falls under the category of interactive imitation learning. By incorporating this approach, the learned agent is capable of scaling up and achieving improved performance on larger problem instances. 

In recent works, 
\citet{yilmaz2021study} also trained a node selection policy with imitation learning. However, their approach uniquely introduced a node comparison operator, which determines whether to expand the left child, the right child, or both children of a given node. The operator can be combined with a backtracking algorithm to provide a full node selection policy.
\citet{labassi2022learning} combined the imitation learning framework introduced by \citet{He2014} with the bipartite graph representation of mixed-integer programming problems developed by \citet{Gasse2019}. Their method showed improved performance compared to the previous methods.

\begin{table}
    \centering 
        \begin{tabular}{p{0.13\linewidth} | p{0.25\linewidth} | p{0.25\linewidth} | p{0.28\linewidth} }
        \toprule
        Paper & Goal  & Methodology & Features\\
        \midrule
        \cite{He2014} &  Learn to mimic a theoretically optimal node selection strategy which knows the optimal solution in advance and only expands nodes containing the optimal solution & Imitation learning via SVM & node features, branching features and B\&B tree features\\
        \midrule 
        \cite{song2018learning} &  Same as \cite{He2014} & Retrospective imitation learning & node features and tree features \\
        \midrule
        \cite{yilmaz2021study} & Learn to mimic the SCIP's default node selection strategy & Use imitation learning to learn an operator that can decide which child node to expand & Same as \cite{Gasse2019}\\
        \midrule
         \cite{labassi2022learning} & Same as \cite{He2014} & Imitation learning via GNN & Same as \cite{Gasse2019} \\
        \bottomrule 
        \end{tabular}
    \caption{Summary of works using AI techniques to enhance node selection in the B\&B algorithm.} \label{tab:b&b_ns} 
\end{table}

\subsection{Mixed integer non-linear programming} \label{sec:minlp}
The field of mixed-integer nonlinear programming (MINLP) has seen some progress in the integration of AI techniques, however, it remains an area with less maturity compared to mixed-integer linear programming. One example of the use of AI in MINLP is the work of \citet{baltean2019scoring}, who attempted to learn linear outer approximations of semidefinite constraints for non-convex quadratic programming problems with box constraints. In this study, a neural network was used to select the most promising submatrices without incurring the computational overhead of solving semidefinite programming problems to generate the necessary cuts. \citet{ghaddar2022learning} focused on using AI to select the ``best branching strategy" in the context of a B\&B search tree embedded within the reformulation-linearization technique (RLT) for solving polynomial problems. They designed hand-crafted features to select the branching strategy that optimizes a quantile regression forest-based approximation of their performance indicator. Additionally, \citet{gonzalez2022polynomial} considered a portfolio of second-order cone and SDP constraints to strengthen the RLT formulation for polynomial problems and used AI to select constraints to add within a B\&B framework.

In a related line of research, several studies have explored the use of AI to predict the computational advantages of certain techniques for solving mixed-integer quadratic programs (MIQPs) and nonconvex mixed-integer nonlinear programs. \citet{bonami2018learning} trained classifiers to predict the computational benefits of linearizing products of binary variables or binary and bounded continuous variables for solving MIQPs. \citet{nannicini2011probing} used support vector machine (SVM) classification to decide whether an expensive optimality-based bound tightening routine should be used instead of a cheaper feasibility-based routine for nonconvex MINLPs.

\subsection{Enhancement for cutting-plane methods} \label{sec:cutting_plane}

Cutting-plane methods is one of the famous methods for solving integer programming problems \citep{gomory1960algorithm}. Consider the following integer programming problem
\begin{alignat*}{6}
    \min_{x\in\mZ^n} \enspace & \enspace c^Tx \\
    \text{s.t.} \enspace & \enspace Ax \leq b \\
                         & \enspace x \geq 0.
\end{alignat*}
Let 
$\Cscr = \{x \in \mZ^n \mid Ax \leq b, \enspace x \geq 0\}$ denote the feasible region for the integer programming problem. 
The cutting plane method starts with solving the LP obtained from the above problem by dropping the integrality constraints $x \in \mZ^n$. Let 
$\Cscr^{(0)} = \{x \in \mR^n \mid Ax \leq b, \enspace x \geq 0\} \supseteq \Cscr$ 
denote the feasible region for the relaxed linear programming problem. Let $x^{(0)} \in \Cscr^{(0)}$ denote the optimal solution to the relaxed LP problem. Let's assume $x^{(0)} \notin \mZ^n$. 
The cutting plane method then finds a cut $(\alpha^{(0)}, \beta^{(0)})$ such that 
\[
    {\alpha^{(0)}}^T x \leq \beta^{(0)} \enspace \forall x \in \Cscr
    \enspace \text{and} \enspace
     {\alpha^{(0)}}^T x^{(0)} > \beta^{(0)}
\]
Then the new feasible region is constructed by 
\[\Cscr^{(0)} \supseteq \Cscr^{(1)} = \Cscr^{(0)} \cap \{x\in\mR^n \mid {\alpha^{(0)}}^T x \leq \beta^{(0)}\} \supseteq \Cscr,\]
and the corresponding LP is solved to obtain $x^{(1)}$. This procedure iterates until $x^{(t)} \in \mZ^n$, which can be shown to be the optimal solution for the original integer programming problem. 

Gomory cuts are a typical way of generating cuts \citep{gomory1960algorithm}. Let $x^{(t)}$ denote the current iterate and define 
\[I^{(t)} = \{i \in [n] \mid x^{(t)}_i \notin \mZ\}.\]
Then for each $i \in I^{(t)}$, we can generate a cut by having
\[ \alpha_i^{(t)} = -A^{(t)}_i + \lfloor A^{(t)}_i \rfloor 
\enspace\text{and}\enspace
\beta_i^{(t)} = -b^{(t)}_i + \lfloor b^{(t)}_i \rfloor \]
where $A^{(t)}$ and $b^{(t)}$ are the constraints in the current LP. Then we choose one of the possible Gomory cuts 
\[\Dscr^{(t)} = \{(\alpha_i^{(t)}, \beta_i^{(t)}) \mid i \in I^{(t)}\}\]
to add to the current LP. The selection of the Gomory cut largely affects the performance of the cutting-plane method. \citet{tang2020reinforcement} propose a reinforcement-learning-based approach to select cuts. At iteration $t$, the state space is given by the current LP and the corresponding optimal solution, i.e., 
\[s^t = (\Cscr^{(t)}, c, x^{(t)}).\]
The available actions are given by 
\[a^t = \Dscr^{(t)},\] 
consisting of all possible Gomory's cutting planes, and the reward at time step $t$ is given by the increase in the objective value, i.e., 
\[R^t = c^Tx^{(t+1)} - c^Tx^{(t)}.\]
Finally, we describe the design of the policy network $\pi_\theta(a^t \mid s^t)$. In order to handle problem instances of different sizes, the authors utilize the LSTM network architecture. More specifically, they first embed all the constraints in the current LP and all the candidate cuts into the same space via an LSTM network, 
\begin{equation*}
    \begin{split}
        h_i &= LSTM_\theta([a_i, b_i]), \enspace \forall [a_i, b_i] \in \Cscr^{(t)} \\
        g_j &= LSTM_\theta([\alpha_j, \beta_j]), \enspace \forall [\alpha_j, \beta_j] \in \Dscr^{(t)}.
    \end{split}
\end{equation*}
Then the score for every candidate $[\alpha_j, \beta_j] \in \Dscr^{(t)}$ is computed by 
\[\text{score}_j = \sum_{[a_i, b_i] \in \Cscr^{(t)}} h_i^T g_j.\]
Finally, $\pi_\theta(a^t \mid s^t)$ returns the probabilities over the action space by applying a softmax function to all the computed scores. 

Similarly, \citet{huang2022learning} also propose a cut-selection strategy based on multiple-instance learning. The training process is similar to reinforcement learning, and the major difference is that the reward is defined by the reduction in the total running time. Alternatively, \citet{paulus2022learning} design another cut-selection strategy via imitation learning. The key idea is that they design a greedy selection rule which is called {looking ahead}. More specifically, for each candidate cut 
\[[\alpha, \beta] \in \Dscr^{(t)},\]
the cut is added to the current LP and solve for the optimal solution 
\[x^{(t)}_{\alpha, \beta} = \argmin\{c^T x \mid x \in \Cscr^{{t}}, \ \alpha^T x \leq \beta\},\]
and the looking ahead rule will select the cut that improves the objective most, i.e.,
\[[\alpha^{(t)}, \beta^{(t)}] \in \argmax\{s_{\alpha, \beta} \coloneqq c^T x^{(t)}_{\alpha, \beta} - c^T x^{(t)} \mid [\alpha, \beta] \in \Dscr^{(t)}\}.\]
Looking ahead is a strong rule for selecting cuts—but an expensive one. At every iteration, running it requires to solve $|\Dscr^{(t)}|$ additional LPs. Therefore, the authors propose to use looking ahead scores to facilitate the training of a policy for cut selection via imitation learning. They represent the current LP, and the candidate cuts as a tripartite graph whose nodes are divided into three parts: variables, constraints and cuts. Then they use the standard GNN approach to predict the looking ahead score for each cut node with a soft binary entropy loss.

\subsection{Heuristics} \label{sec:heuristics}


For discrete optimization problems, besides exact methods like Branch-and-Bound and Branch-and-Cut, heuristic algorithms are also widely adopted for their simplicity and efficiency. Although they do not guarantee to find the optimal solution, they are designed to find a good solution within reasonable computational time.
Traditional heuristics can be classified into two categories based on their goals: 1) Finding a {feasible} solution quickly. The feasibility pump (FP) \citep{berthold2019ten} and diving \citep{Maniezzo2021} heuristic fall into this category, 2) Trying to find a {high-quality} feasible solution, even though there's no assurance of achieving optimality.
Local branching~\citep{Fischetti2003LocalB}, relaxation induced neighborhood search (RINS)~\citep{Danna2005ExploringRI}, and large neighborhood search (LNS)~\citep{Shaw1998UsingCP,Pisinger2018LargeNS} fall into this category. 
It's essential to note that these two categories can be combined, i.e., find an initial feasible solution quickly and then enhance it. For instance, one might initially employ FP and then RINS to refine the solution further.


AI-based heuristics often draw inspiration from traditional heuristics. 
In the following, we first discuss the heuristics inspired by category 1) and then by category 2).
Within category 1), two AI-based heuristics have been proposed inspired by FP and diving respectively. 
\citet{qi2021smart} propose a smart feasibility pump (SFP) method inspired by the traditional feasibility pump (FP) heuristic. This new method surpasses the traditional one by reducing the number of steps needed to reach the first feasible solution. FP first solves the relaxed problem~\eqref{eq:mip-relax}. Its solution ${x}^{(0)}$ is then rounded to obtain an initial integer solution $\bar{x}^{(0)}$. However, due to rounding, this solution is usually infeasible w.r.t. constraints~\eqref{eq:mip-b} and \eqref{eq:mip-c}. 
FP iteratively projects the solution $\bar{x}^{(t)}$ back to the feasible region of the relaxed problem~\eqref{eq:mip-relax}, obtaining $x^{(t+1)}$ and rounding it to $\bar{x}^{(t+1)}$. 
The projection wants to satisfy constraints~\eqref{eq:mip-b} to \eqref{eq:mip-c} with the smallest change on the current solution, while the rounding focuses only on the integer constraints. 
This process may take a long time since the projection and rounding focus on short-term benefits. Furthermore, there is a risk of stagnation, i.e., failing to find any feasible solution. 
In contrast, SFP trains an RL agent such that the infeasibility will rapidly decrease in several iterations. 
SFP inherits the framework of FP, but the key point is that SFP allows the agent to change the current solution $x^{(t)}$ before rounding. 
The {reward} is defined as negative infeasibility. The cumulative rewards will drive the agent to balance between current-step and long-term negative infeasibility. 
The {state} includes the information of Problem~\eqref{eq:mip} and the current solution. 
The {action} of the agent is the change $a^{(t)}$ in the current solution $x^{(t)}$. 
Specifically, the {state transition} incorporates the change and rounding, {i.e.}, $x^{(t+1)} = [x^{(t)} + a^{(t)}]$, where $[\cdot]$ is the rounding operation. 
Despite the non-differentiable nature of rounding and infeasibility with respect to the agent's actions, the PPO algorithm (See section \ref{sec:rl}) can still be effectively utilized, since it is designed to handle such a challenge. 

Alternatively, \citet{nair2021solving} propose neural diving. In the literature, diving is a heuristic that explores the branch-and-bound tree.
It starts from any partial assignment of integer variables and finds a feasible assignment of integer variables (if it exists) in an iterative way. 
In contrast, neural diving wants to directly {predict} an initial assignment of integer variables. 
This assignment may not be feasible, but needs to be high-quality, i.e., resulting in a high objective value. This distinguishes neural diving from traditional diving techniques since the goal has changed. 
The prediction of the assignment is achieved by a GNN inference. 
The advantage is that this inference is a ``one-off deal''.
Next, we discuss how the GNN model is trained. 
Let $P_i$ represent the $i$-th MIP problem in the training set. Each problem has its graph representation. 
If we solve the MIP with the B\&B algorithm, {multiple} feasible solutions can be found on the leave nodes, and each solution corresponds to an assignment of integer variables. 
We denote the number of these feasible assignments as $N_i$ and the set of these feasible assignments as $\Xscr_i=\{{x_{(i,j)}}\}_{j=1}^{N_i}$. Here, $x_{(i,j)}$ is the $j$-th assignment of integer variables for $i$-th MIP problem. 
Now the quality of $j$-th assignment of integer variables is evaluated by the following energy score $E\left(x_{(i,j)}; P_i \right)$. 
$$
E\left(x_{(i, j)} ; P_i\right)= \begin{cases}\hat{f}\left(x_{(i, j)}\right) & \text { if } x_{(i, j)} \text { is feasible } \\ \infty & \text { otherwise }\end{cases}
$$
where $\hat{f}(x_{(i,j)})$ 
represents the objective value obtained by assigning $x_{(i,j)}$ to integer variables in $P_i$ and assigning the rest continuous variables according to the solution of the resulting linear program. 
This energy score is normalized on all assignments for $i$-th problem and converted to sample weight $w_{ij}$. 
To summarize, now we have a training set $\{(P_i, \Xscr_i)\}_{i=1}^{N}$. We can train the neural network $g(\cdot;\theta)$ with parameter $\theta$ with the following training loss.  
$$
\Lscr(\theta)=-\sum_{i=1}^N \sum_{j=1}^{N_i} w_{i j} \ell \left( g( P_i; \theta), x_{(i, j)} \right) ,
$$
where the loss function $\ell(\cdot,\cdot)$  measures the difference between the prediction and the ground truth. 
The prediction $g( P_i; \theta)$ is of high quality since it is highly likely to be feasible and close to optimal. 
It is guaranteed that the provided $x_{(i,j)}$ in {training set} are all feasible assignments. However, 
$g( P_i; \theta)$ may still provide an infeasible assignment at the inference stage.  
Thus, neural diving is usually followed by a search method that allows infeasibility, such as B\&B \citep{nair2021solving} and large neighborhood search \citep{song2018learning, sonnerat2021learning}.

For category 2), traditional heuristics are typically iterative, where early decisions can significantly influence long-term outcomes. Thus, reinforcement learning is frequently employed to enhance these traditional heuristics.
However, it's important to note that RL, while powerful, still encounters convergence issues, especially when confronted with large action and state spaces, or sparse reward signals. 
A prevalent approach to mitigate these challenges involves initially employing imitation learning to train the agent, followed by fine-tuning using reinforcement learning. 


\citet{Khalil2017} focus on challenging combinatorial optimization problems on graphs, e.g., the Minimal Vertex Cover (MVC) problem. This problem seeks the smallest subset of vertices such that each edge in the graph is incident to at least one vertex in this subset. 
To solve this problem, traditional greedy heuristics progressively build the desired subset by incorporating the vertex perceived as the most advantageous during each iteration.
\citet{Khalil2017} employ reinforcement learning to design AI-enhanced greedy heuristics. 
The ``action'' refers to the selection of a vertex into the subset. The action space includes the vertices that satisfy the graph problem's constraints. 
The ``state" contains the current subset and the graph embedding. Here, ``graph embedding" is a numerical vector that summarizes the information in the graph. 
The ``reward" is directly the objective function of the original problem. 
In this way, the learned agent doesn't merely seek a feasible solution but aims for the solution with high objective value. 
Meanwhile, the agent is equipped to handle ``delayed rewards'' (See Section \ref{sec:rl}), which are challenging for traditional greedy algorithms. 

In the following, we introduce two AI-enhanced heuristics \citep{sonnerat2021learning,song2018learning}, building upon existing Large Neighborhood Search (LNS) heuristics. We start by briefly introducing LNS.  
LNS is a combinatorial optimization heuristic that begins with an assignment of values for variables and iteratively refines it by searching a large neighborhood around the current assignment. 
In each iteration, LNS evaluates multiple neighborhood solutions, comparing which one is most promising for rapidly converging to a high-quality solution. 
This naturally raises two questions: 
1) How to find neighborhood solutions?
2) How to assess whether a particular neighborhood is more promising than another? For question 1), existing heuristics have a consistent approach. 
From the current solution, they ``unassign" specific variable values, marking them as ``unknown." The resulting MIP is then exactly solved using off-the-shelf solvers to obtain the neighborhood solution.
AI-assisted heuristics inherit this approach but will smartly decide which variable values to ``unassign". For question 2), a naive method leverages local information. Specifically, neighborhood solutions are evaluated using metrics like objective value or primal gap. The most promising neighborhood is selected. 
This method, however, has two main limitations: a) computational cost and b) its short-sighted nature.
Limitation a) exists because each neighborhood solution is obtained by solving an MIP. 
To overcome this limitation, imitation learning is utilized. 
The benefit is that during the inference stage, the learned AI model directly predicts the promising neighborhood, which is a ``one-off'' deal.  
Regarding limitation b), the term ``short-sighted" refers to the tendency to make ``greedy" decisions based purely on local information. Simply identifying the locally best neighborhood does not guarantee convergence to a high-quality solution in the long run.
RL provides a remedy to this challenge due to deleted rewards. 
In the end, 
the common paradigm is to initially utilize imitation learning to train an agent, followed by refinement through RL. 
The purpose of this paradigm is not merely to address the aforementioned limitations but also to help the learning process of agents converge. 
Since after imitation learning the agent is already able to make a rational decision, in later RL stages, the agent needs less exploration and thus converges faster. 

Two AI-enhanced LNS heuristics follow the same paradigm, with the primary distinction being their design of action. \citet{sonnerat2021learning} design the action directly as which integer variables to unassign while
\citet{song2020general} design the action as a decomposition. 
Specifically, it means decomposing the integer variable sets \(x\) into disjoint subsets, i.e., \( x = x_1 \cup x_2 \cup ... \cup x_k ... \cup x_K \). Here, the hyper-parameter \(K\) represents the number of equally sized subsets. 
This decomposition action is predicted by an AI model (multi-layer perceptron). The model predicts the subset index to which each variable should be allocated. 
As for state transition, a number of \(K\) different MIPs are solved. 
For the \(k\)-th MIP (\(k \in \{1,...,K\}\)), variables outside of the subset \(x_k\) retain their values from the current solution, while an optimization solver refines the variables within \(x_k\). 
The solution with the best primal gap is then selected. A shared limitation in both AI-enhanced heuristics is the lack of adaptive control over the neighborhood size. In the work by \citet{sonnerat2021learning}, the neighborhood size is determined by a hyperparameter representing the number of unassigned integer variables. In \citet{song2020general}, the neighborhood size is controlled by the number of subsets \(K\). They are both static hyper-parameters.
Although larger neighborhoods can reduce the risk of trapping to local optima, this also increases computational cost. A promising direction for future research is the development of mechanisms for adaptive neighborhood sizing.

To conclude, AI models assist or surpass heuristics either by imitating a computationally expensive target or by learning from the costly interaction with the environment. At the inference stage, AI models can quickly take smart actions and even outperform the state-of-the-art commercial solvers on MIP with hundreds or thousands of variables \citep{song2020general}.

\section{Conclusion} \label{sec:conclusion}
AI for Operations Research typically focuses on enhancing an individual stage within the Operations Research pipeline. However, exploring interactions between different stages is an intriguing area of study, such as gathering feedback from later stages to refine earlier ones. The smart predict-then-optimize paradigm is one example as a pioneering approach for integrating diverse stages. 
Under the SPO framework, besides minimizing the objective, an interesting extension is also ensuring feasibility. Relaxing feasibility is a widely recognized technique employed in real-world applications. Its motivation stems from the fact that domain experts often find the response of ``impossible to complete the task" unsatisfactory. Instead, they strive to comprehend the reasons for infeasibility and explore avenues to relax constraints slightly, thereby achieving a feasible solution. 
By incorporating this concept into SPO, it is interesting to let the AI model be capable of predicting parameters or adjusting parameters such that the constraint is feasible. 
Furthermore, other types of interactions can also prove valuable. For instance, if the optimization process exhibits slow convergence or sub-optimality, it may indicate the need for an alternative formulation. It is beneficial to explore the potential of AI models to automatically perform such adaptations.
Nevertheless, designing a space that encompasses well-known optimization formulations and provides the building blocks for modifying formulations is a challenging task. Achieving a framework that automatically adapts formulations also remains a promising area of research.

Despite the advancements in automatic algorithm configuration and algorithm selection, the task of ``unified software selection \& tuning'' remains interesting and challenging. 
Given an optimization problem instance, and a set of available optimization software, e.g. Gurobi, CPLEX, and OptVerse, we want to predict which software along with its potential configuration is going to perform best on this instance. On the one hand, the challenge is inherited from the algorithm selection task, i.e., the target (algorithm or software) is regarded as a black box. On the other hand, the software incorporates several parameters and preprocessing algorithms that result in a complex structure. To be specific, if CPLEX performs worse than Gurobi on an instance, the root cause might be an inappropriate configuration rather than the software's inadequacy. In other words, the software selection task combines the algorithm selection and automatic algorithm configuration tasks.

In the model-based method for automatic algorithm configuration, a performance model is built for {Algorithm Runtime Prediction}. Specifically, the goal is to employ AI techniques to predict the runtime of an algorithm on a previously unseen problem given the problem-specific features. While the existing work usually claim the performance model is empirically strong, it will be useful to know its boundary, i.e., when it will fail. To deal with this, we may need a measurement for smoothness and learnability of the performance function. We may also need a similarity metric between the unseen input and training instances. 

AI-based methodologies implicitly assume that the testing instances are similar to the ones used in training. However, there currently does not exist a well-established {similarity metric} exists. 
The challenge is that the property of a similarity metric is hard to describe mathematically. We may define the metric by looking at the model formulation, structure, and the corresponding parameters, but it is hard to derive a guarantee that similar problems will benefit similarly from an AI model. 
Nevertheless, designing an empirical similarity metric with AI techniques might still be valuable. 

In conclusion, AI techniques have demonstrated great potential in enhancing each stage of the OR process. More future works are worth exploring the synergy between AI and OR, e.g., utilizing AI to enhance the interactions between different OR stages. 
This synergy will undoubtedly lead to exciting advancements and novel solution methods in a multitude of domains. 


\bibliographystyle{apalike}
\bibliography{biblio}
\end{document}